\documentclass[a4paper,fleqn]{cas-sc}
\usepackage{xcolor}
\usepackage{todonotes}
\usepackage{algorithm}
\usepackage{algpseudocode}
\usepackage{mathtools}
\usepackage{xspace}
\usepackage{subcaption}
\usepackage{standalone}
\usepackage{amsfonts}
\usepackage{amssymb}
\usepackage{amsthm}
\usepackage{pifont}
\usepackage{cleveref}
\usepackage{thm-restate}
\newcommand{\arcs}{\mathcal{A}}
\newcommand{\nodes}{\mathcal{V}}

\newcommand{\dleq}{\leq_{D}}
\newcommand{\dom}{<_D}

\newcommand{\bestreal}[1]{\mathbf{c}[#1]}
\newcommand{\paths}{\mathcal{P}}
\newcommand{\shpaths}{\mathbf{P}}
\newcommand{\pathst}{\paths_{s,t}}
\newcommand{\shpathst}{\shpaths_{s,t}}

\newcommand{\ptn}{\ensuremath{\mathrm{\textit{PTN}}}}

\newcommand{\eshpathst}{\hat{\shpaths}_{s,t}}

\newcommand{\effpaths}{\mathbb{P}_{st}^{\dleq}}

\newcommand{\weffpaths}{\mathbb{P}_{st}^{\dom}}
\newcommand{\paretoset}{\mathbf{W}}

\newcommand{\bl}{\ensuremath{\boldsymbol{\ell}}}
\newcommand{\bu}{\ensuremath{\boldsymbol{u}}}

\newcommand{\dep}{\text{dep}}
\newcommand{\arr}{\text{arr}}

\newcommand{\cispp}{\textnormal{C-ISPP}\xspace}
\newcommand{\eispp}{\textnormal{E-ISPP}\xspace}

\newcommand{\mosp}{\textnormal{MOSP}\xspace}

\newcommand{\mc}{\multicolumn}
\newcommand{\dr}{\textcolor{darkred}}

\newcommand{\timpasslib}{\texttt{TimPassLib}}

\newcommand{\extremepoints}{\ensuremath{\mathcal{E}([\bl, \bu])}}

\newcommand{\Wup}[1]{\texttt{W#1}}
\newcommand{\Cells}{\ensuremath{\mathcal C}}
\newcommand{\Stations}{\ensuremath{\mathcal S}}
\newcommand{\Edges}{\ensuremath{\mathit E}}

\newcommand\mscriptsize[1]{\mbox{\scriptsize\ensuremath{#1}}}

\newcommand{\exn}[1]{\mscriptsize{\times 10^{#1}}}

\newcommand{\cmark}{\ding{51}}%
\newcommand{\xmark}{\ding{55}}%
\definecolor{tabblue}{HTML}{1f77b4}
\definecolor{taborange}{HTML}{ff7f0e}
\definecolor{tabgreen}{HTML}{2ca02c}
\definecolor{tabred}{HTML}{d62728}
\definecolor{tabpurple}{HTML}{9467bd}
\definecolor{dgreen}{rgb}{0.0, 0.26, 0.15}
\definecolor{tabyellow}{rgb}{0.99,0.78,0.07}
\definecolor{darkred}{HTML}{800000}

\usepackage[longnamesfirst, numbers]{natbib}

\def\tsc#1{\csdef{#1}{\textsc{\lowercase{#1}}\xspace}}
\tsc{WGM}
\tsc{QE}

\newtheorem{theorem}{Theorem}
\newtheorem{lemma}[theorem]{Lemma}
\newtheorem{remark}[theorem]{Remark}
\newtheorem{corollary}[theorem]{Corollary}
\newtheorem{definition}[theorem]{Definition}
\newtheorem{example}[theorem]{Example}

\begin{document}
\let\WriteBookmarks\relax
\def\floatpagepagefraction{1}
\def\textpagefraction{.001}

\shorttitle{Computing All Shortest Passenger Routes with a Tropical Dijkstra Algorithm}

\shortauthors{B. Masing, N. Lindner, E. Bortoletto}  

\title [mode = title]{Computing All Shortest Passenger Routes with a Tropical Dijkstra Algorithm}  



%

\author[1]{Berenike Masing}[orcid = 0000-0001-7201-2412] 

\cormark[1]

\fnmark[1,2]

\ead{masing@zib.de}


\credit{}

\affiliation[1]{organization={Zuse Institute Berlin},
            addressline={Takustr. 7}, 
            postcode={14195}, 
            state={Berlin},
            country={Germany}}

\author[1]{Niels Lindner}[orcid = 0000-0002-8337-4387] 


\ead{lindner@zib.de}


\credit{}
\author[1]{Enrico Bortoletto}[orcid = 0000-0002-2869-6498]
\ead{bortoletto@zib.de}
\fnmark[2]

\cortext[1]{Corresponding author}

\fntext[1]{Funded within the Research Campus MODAL, funded by the German Federal Ministry of Education and Research (BMBF) (fund number 05M20ZBM).} 

\fntext[2]{Funded by the Deutsche Forschungsgemeinschaft (DFG, German Research Foundation) under Germany's Excellence Strategy – The Berlin Mathematics Research Center MATH+ (EXC-2046/1, project ID: 390685689).}


\begin{abstract}
Given a public transportation network, which and how many passenger routes can potentially be shortest paths, when all possible timetables are taken into account? 
This question leads to shortest path problems on graphs with interval costs on their arcs and is closely linked to multi-objective optimization.
We introduce a Dijkstra algorithm based on polynomials over the tropical semiring that computes complete or minimal sets of efficient paths.
We demonstrate that this approach is computationally feasible by employing it on the public transport network of the city of Wuppertal and instances of the benchmarking set \timpasslib, and we evaluate the resulting sets of passenger routes.
\end{abstract}




\maketitle


\section{Introduction}
\label{sec:intro}

In a public transportation network, passengers typically have several options to choose routes for their journeys. The route choice is heavily influenced by the timetable \cite{schmidt_integrating_2014}. However, the intuition is that only a small fraction of all paths in the network will be reasonable passenger routes, regardless of the timetable. If this can be confirmed, then this opens up new algorithmic perspectives for integrated timetabling and passenger routing \cite{schiewe_introducing_2023}, where one bottleneck is the perpetual computation of shortest paths for a large number of origin-destination (OD) pairs \cite{lobel_restricted_2020,schiewe_periodic_2020}. 

In this paper, we want to corroborate this intuition in the following setting: We will consider directed graphs whose arcs are labelled with real intervals. This way, we can consider event-activity networks as used in periodic timetabling, where the duration of the activities is still unknown, but general parameters such as minimum driving times or minimum transfer times have already been determined. On such a graph, we ask for the set of all paths that are a shortest path for at least one scenario of costs within the interval bounds, and call this the Complete Interval Shortest Path Problem (\cispp). This set of paths is the set of all potential shortest passenger routes on the level of event-activity networks. For the application of integrated timetabling and passenger routing, it is however sufficient to know only one shortest path for each cost scenario. We therefore aim to also compute a minimal scenario-covering set of paths in the Essential Interval Shortest Path Problem (\eispp). If this set turns out to be rather small, an enumeration could be used to decrease the complexity of current methods.

Both \cispp and \eispp are closely related to robust shortest path problems, where also graphs with unknown costs are considered \cite{catanzaro_reduction_2011,karasan_robust_2012,montemanni_exact_2004}, and multi-objective shortest path problems, as we ultimately complete sets of efficient paths \cite{ehrgott,maristany_de_las_casas_new_2024}, but with respect to an objective whose dimension is exponential in the size of arcs. A preprocessing approach for graphs with interval costs has been presented in \cite{hamonic_optimizing_2023,lindner_optimal_2021}.

To solve \cispp and \eispp efficiently, we propose a novel tropical Dijkstra algorithm that uses polynomials over the tropical min-plus semiring \cite{tropicalGeom_sturmfels,joswig_essentials_trop} as labels. This has the advantage that dominance checks are mere evaluations of tropical polynomials. A different, but also tropical geometry-inspired breadth-first search technique has already been applied to enumerate parametric shortest path trees, but only with few variable arcs \cite{JoswigSchroeter_2022}.

Our tropical Dijkstra algorithm applies to any graph with interval costs. We evaluate our implementation of the method on a set of realistic instances based on the network of the city of Wuppertal, as well as instances derived from \timpasslib, a benchmarking library for the integrated periodic timetabling and passenger routing problem \cite{schiewe_introducing_2023}. It turns out that it is computationally feasible to solve \cispp and \eispp even on the largest instances.
We finally present some statistics on the number of passenger routes.

\Cref{sec:ispp} introduces \cispp and \eispp, and develops their general theory. The connection to multi-objective shortest path problems is explored in \Cref{sec:mosp}. 
Tropical polynomials and the tropical Dijkstra algorithm are presented in \Cref{sec:tropical}. 
For our application of shortest passenger routes in public transport, we illustrate our modelling of event-activity networks in \Cref{sec:application}.
We evaluate our tests in \Cref{sec:experiments}, before concluding the paper in \Cref{sec:conclusion}. Any missing proofs can be found in the appendix. 

\subsection{Our Contributions}
Motivated by shortest passenger paths for not-yet-determined timetables, we formally define both the complete and essential set of shortest paths for a graph with interval costs. We show the equivalence between finding the complete and essential set of shortest paths to a multi-objective shortest path problems with $2^{|\arcs|}$-dimensional objective functions. Using a well-established description of parametric shortest paths as tropical polynomials, we develop an algorithm, Tropical Dijkstra, which can compute both the complete and essential set of shortest paths. The algorithm proves to be effective in computational experiments on a diverse set of realistic instances. Our results show the importance of integrating passenger routing with timetabling and provide the means to pre-compute relevant passenger paths. 


\section{The Interval Shortest Path Problem}
\label{sec:ispp}

Throughout this paper, we will consider the following setup:
Let $G$ be a directed graph with node set $\nodes(G)$ and arc set $\arcs(G)$, and let $\ell, u \in \mathbb{R}^{\arcs(G)}$ be bounds on the arcs with $0 \leq \ell_a \leq u_a$.
The arcs $a \in \arcs(G)$ are hence labelled with intervals $[\ell_a, u_a]$, and we will therefore call $(G, \ell, u)$ a \emph{graph with interval costs}.

Moreover, we call $c \in [\bl, \bu] \coloneqq \times_{a\in \arcs(G)} [\ell_a, u_a]$ a \emph{scenario}, and by $\extremepoints$ we denote the extreme points of the scenarios, i.e., all points of the form $(x_a)_{a\in\arcs(G)}$ with $x_a \in \{\ell_a, u_a\}$.  
A \emph{path} $p$ in $G$ is an alternating sequence $(v_0,a_1,v_1,\ldots, a_l,v_l)$ of nodes $v_i \in \nodes(G)$ and arcs $a_i =(v_{i-1}, v_i) \in \arcs(G)$.
When speaking of paths we will always assume them to be simple, we will denote by $\arcs(p)$ the set of arcs of the path $p$, and will -- for better readability -- describe a path by its sequence of arcs only. 
We denote by $\pathst(G)$ the set of all $s$-$t$-paths in $G$, for $s, t \in \nodes(G)$, and we write $\pathst$ if the graph is clear from context.
For $p \in \pathst$ and a fixed scenario $c \in [\bl,\bu]$, we define the cost of $p$ as $c(p) \coloneqq \sum_{a\in \arcs(p)} c_a$. 
A path $p \in \pathst$ is a \emph{shortest path} of $(G, \ell, u)$ if it is a path of minimum cost w.r.t.~some scenario $c\in [\bl, \bu]$, i.e., $c(p) \leq c(p')$ for all $p'\in \pathst$. 

In the following, we will discuss two flavours of the shortest path problem in graphs with interval costs, namely one in which we are interested in a complete set of shortest paths, and one in which we want to find only an essential set of shortest paths.

\subsection{The Complete Interval Shortest Path Problem}
\label{subsec:cispp}

Let $(G, \ell, u)$ be a graph with interval costs. For $s, t \in \nodes(G)$ and $c \in [\bl, \bu]$, we denote by $\shpathst^c$ the set of all shortest $s$-$t$-paths in $G$ w.r.t.\ $c$.

\begin{definition}[Complete Interval Shortest Path Problem (\cispp)]
    For a graph with interval costs $(G,\ell, u)$ and a source node $s\in \nodes(G)$, the \emph{Complete Interval Shortest Path Problem (\cispp)} consists in finding the collection of $s$-$t$-paths $\shpathst \coloneqq \bigcup_{c \in [\bl, \bu]} \shpathst^c\subseteq \pathst$, for all $t\in \nodes(G)$.
\end{definition}

\begin{figure}
    \centering
     \begin{tikzpicture}[xscale=3, yscale = .7]
		\tikzstyle{a} = [line width=0.8, ->, black]
            \tikzstyle{n} = [draw, circle, inner sep = 1.5, minimum width = 10, black]
            \node[n] (s) at (0,0) {$s$};
            \node[n] (v) at (1, 0) {$v$};
            \node[n] (t) at (2,0) {$t$};
            \draw[a] (s) to[bend left=60, looseness=1.2] node[above, midway]{$a_1$ : [5,6]} (v);
            \draw[a] (v) to[bend left=60, looseness=1.2] node[above, midway]{$a_2$ : [4,6]} (t);
            \draw[a] (s) to[bend right=60, looseness=1.2] node[below,midway] {$a_3$ : [2,6]} (v);
            \draw[a] (v) to[bend right = 60, looseness = 1.2] node[below,midway] {$a_4$ : [2,6]} (t);
            \draw[a] (s) to[bend right = 90, looseness = 1.8] node[below,midway] {$a_5$ : [5,8]} (t);
            \draw[a] (s) to[bend left = 90, looseness = 1.8] node[above,midway] {$a_6$ : [8,12]} (t);
 \end{tikzpicture}
 
    \caption{An example of an interval shortest path instance $(G, \ell, u, s)$, with labels $a \colon [\ell_a, u_a]$.}
    \label{fig:running}
\end{figure}
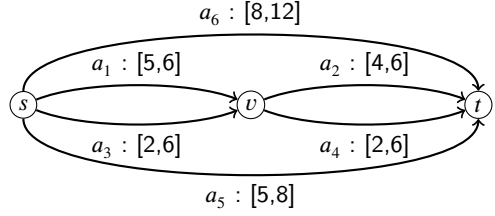

\begin{example}
    Consider the graph depicted in \Cref{fig:running}.
    This instance has six $s$-$t$-paths, all of which are shortest paths for some scenario, except for the path $(a_1, a_2)$, which at best has cost $9$, and is therefore always more expensive than taking $a_5$, which has cost of at most 8. 
    See \Cref{tab:path_table_of_running} with the list of all $s$-$t$-paths and a corresponding scenario certificate for each path that is in $\shpathst$. The remaining information in the table will be explained later on. 
\end{example}

\begin{table}
    \centering
    \begin{tabular}{l | cl c l}
        path $p$ & $\in \shpathst$ &   scenario certificate $\bestreal{p}$ & $\in \eshpathst$ & tropical path monomial $\tau(p)$\\
        \hline
        $p_1 = (a_1, a_2)$ & \xmark&   -- & \xmark & $x_{1} \odot x_{2}$\\
        $p_2 = (a_1, a_4)$ & \cmark &  $(5,6,6,2,8,12)$& \cmark & $x_{1} \odot x_{4}$ \\
        $p_3 = (a_3, a_2)$ & \cmark & $(6,4,2,6,8,12)$& \cmark & $x_{3} \odot x_{2}$\\
        $p_4 = (a_3, a_4)$ & \cmark & $(6,6,2,2,8,12)$& \cmark & $x_{3} \odot x_{2}$\\
        $p_5 = (a_5)$ & \cmark & $(6,6,6,6,5,12)$& \cmark& $x_{5}$\\
        $p_6 = (a_6)$ & \cmark & $(6,6,6,6,8,8)$& \xmark &$x_{6}$
    \end{tabular}
    \caption{All $s$-$t$-paths of \Cref{fig:running}, displaying whether they are in the complete and in an essential set of path, including a corresponding scenario, namely the best scenarios, which certifies them as a shortest path if applicable.}
    \label{tab:path_table_of_running}
\end{table}

\begin{definition}[Best Scenario of a Path]
    Given a path $p \in \pathst$ we say $\bestreal{p} \in \extremepoints$ is the \emph{best scenario} of path $p$, where $\bestreal{p}_a \coloneqq \ell_a$ for $a \in \arcs(p)$ and $\bestreal{p}_a \coloneqq u_a$ for $a\in \arcs(G)\setminus \arcs(p)$.
\end{definition}

The role of the best scenario is underlined by the following technical lemma, that will be the basis for our upcoming results, and whose proof can be found in the appendix.

\begin{restatable}[]{lemma}{technical}
\label{lem:technical}
    Let $p, p' \in \pathst$ and let $\vartriangleleft\, \in \{<, \leq\}$.
    \begin{enumerate}
        \item If $c(p) \vartriangleleft c(p')$ for some $c \in [\bl, \bu]$, then $\bestreal{p}(p) \vartriangleleft \bestreal{p}(p')$.
        \item If $\bestreal{p'}(p) \vartriangleleft \bestreal{p'}(p')$, then $c(p) \vartriangleleft c(p')$ for all $c \in [\bl, \bu]$.
    \end{enumerate}
\end{restatable}

\begin{theorem}
\label{thm:best_scenario_shortest_path}
    A path $p \in \pathst$ is a shortest path for some scenario if and only if it is a shortest path w.r.t.\ its best scenario $\bestreal{p}$. In particular, $\shpathst = \bigcup_{c \in \extremepoints} \shpathst^c$.
\end{theorem}
\begin{proof}
    If $p$ is a shortest path for $c \in [\bl, \bu]$, then \Cref{lem:technical} (1) implies that $p$ is a shortest path for $\bestreal{p}$. It remains to note that $\bestreal{p} \in \extremepoints$.
\end{proof}

The value of the previous theorem is twofold: 
It firstly enables us to verify efficiently whether a given path is a shortest path, simply by choosing the best scenario and running a shortest path algorithm.
Moreover, we can now describe the complete set of shortest paths in terms of the finite set of scenarios $c \in \extremepoints$, as opposed to $c\in [\bl, \bu]$.

\begin{remark}
\label{lem:exponential_shortest_path}
    For each directed graph $G$ there is a choice of $\ell$ and $u$ such that every path is a shortest path in $(G, \ell, u)$: 
    Set, e.g., $\ell_a = 0$ and $u_a =1$ for all $a\in \arcs(G)$. 
    Now consider any $s$-$t$-path $p$ in $G$.
    Clearly, $\bestreal{p}(p) = 0$. Since any other path $p' \neq p$ must use some arc not on $p$, we thus obtain that $\bestreal{p}(p') \geq 1 > 0 = \bestreal{p}(p)$.
    As $p'$ was arbitrarily chosen, this certifies $p$ as a shortest path.
\end{remark}

While \Cref{lem:exponential_shortest_path} shows that the complete set of shortest paths can be exponential in size, it is unlikely in our practical application. To this end, we will proceed to introduce a dominance relation of paths, which will enable us to later describe an algorithm in a concise way.

\begin{definition}[Strict Dominance]
    Let $p, p'$ be two $s$-$t$-paths in $(G, \ell, u)$.
    Then $p$ \emph{strictly dominates} $p'$ -- denoted by $p\prec p'$ -- if and only if $c(p) < c(p')$ for all $c\in [\bl, \bu]$. 
\end{definition}

Using \Cref{lem:technical}, we can easily verify whether one path dominates the other:

\begin{lemma}\label{lem:property_bestreal_dominance}
    Let $p, p'\in \pathst$. 
    Then $p \prec p'$ if and only if $\bestreal{p'}(p) < \bestreal{p'}(p')$.
\end{lemma}

This dominance relation of paths allows us to characterize when a path belongs to $\shpathst$:

\begin{theorem}
\label{thm:equivalence_nonshortestpath}
    Let $p\in \pathst$. 
    The following are equivalent:
    \begin{enumerate}
        \item \label{it:1}
            $p\in \shpathst$.
        \item \label{it:2}
            For all $p'\in \pathst$ holds $p' \not\prec p$ .
        \item \label{it:3}
            For all $p' \in \pathst$ holds $\bestreal{p}(p) \leq \bestreal{p}(p')$. 
    \end{enumerate}
    In particular, \cispp is solved by finding the set of all non-dominated  paths.
\end{theorem}
\begin{proof}
    \eqref{it:2} and \eqref{it:3} are equivalent by \Cref{lem:property_bestreal_dominance}, \eqref{it:1} is equivalent to \eqref{it:3} by \Cref{thm:best_scenario_shortest_path}.
\end{proof}

\begin{remark}
    One could ask the question whether the search for all shortest paths
    is doable, particularly in light of the threat of this set being exponential in size. 
    One idea could be to find a way to compute the graph $G'$ with $\nodes(G') = \nodes(G)$ and $\arcs(G') = \bigcup_{p\in \shpathst} \arcs(p)$, in hopes that all its paths are shortest paths (just like in the standard shortest path problem), or to at least reduce the instance size.
    It is clear that $\shpathst(G) = \shpathst(G')$ when considering the corresponding interval shortest path instance $(G', \ell', u')$ with $\ell'_a = \ell_a$ and $u'_a = u_a$ for all $a\in \arcs(G')$. 
    However, it turns out that not much can be gained from this approach and the graph $G'$: Although there are polynomial-time algorithms that can detect whether an arc is part of a shortest path tree \cite{hamonic_optimizing_2023,lindner_optimal_2021}, it is possible that $G' = G$, as is the case for instance in \Cref{fig:running}. 
    Moreover, there are paths in $G'$ which are not shortest paths, but whose arcs are each part of other shortest paths. 
    In the example, this is the case for the path $(a_1, a_2)$. 
\end{remark}

\subsection{The Essential Interval Shortest Path Problem}
\label{subsec:eispp}

For the application of computing shortest passenger paths in a transportation network for all possible timetable scenarios, a planner might not be interested in the complete set of passenger paths, but would rather have a set of paths which guarantees them to have a shortest path for each timetable realisation. This motivates the following variation of the problem. 

\begin{definition}[Essential Interval Shortest Path Problem (\eispp)]
    For a graph with interval costs $(G,\ell, u)$ and a source node $s\in \nodes(G)$, the \emph{Essential Interval Shortest Path Problem (\eispp)} consists in finding the collection of $s$-$t$-paths $\eshpathst\subseteq \pathst$ such that it is \emph{scenario-covering}, i.e., for each scenario $c\in [\bl, \bu]$ there exists a path $p\in \eshpathst$ which is a shortest path w.r.t.~$c$, and each path in $p\in \eshpathst$ is \emph{essential}, by which we mean that $\eshpathst\setminus\{p\}$ is not scenario-covering.
    We call $\eshpathst$ an \emph{essential set} of shortest paths. 
\end{definition}

It is worth pointing out that an essential set of shortest paths is not necessarily unique. 
For example, if there are two paths $p, p' \in \shpathst$ which have the same costs for all scenarios, then either $p$ or $p'$ can be in an essential set, but not both simultaneously (see also \Cref{thm:equivalence_nonshortest_weak_path}).
In this light, it is natural to then also ask not only for an essential set, but also to attempt to find an essential set of minimum size. We will show that this is automatic in \Cref{sec:mosp}.

We will first develop a few tools to determine whether a path can be in an essential set. 

\begin{definition}[Weak Dominance]
    Let $p, p'\in \pathst$. We say that $p$ \emph{weakly dominates} $p'$ -- denoted by $p \preceq p'$ -- if and only if $c(p) \leq c(p')$ for all $c\in[\bl, \bu]$ and there exists a $c'\in [\bl, \bu]$ with $c'(p) < c'(p')$.  
\end{definition}
We will proceed with mirroring \Cref{lem:property_bestreal_dominance} and \Cref{thm:equivalence_nonshortestpath} for the relation $\preceq$. The first lemma below follows directly from \Cref{lem:technical}, whereas the other proofs can be found in the appendix. 

\begin{lemma}\label{lem:bestreal_weak_dom}
    Let $p, p'\in \pathst$. Then $p \preceq p'$ if and only if $\bestreal{p'}(p) \leq \bestreal{p'}(p')$ and $\bestreal{p}(p) < \bestreal{p}(p')$.
\end{lemma}

\begin{definition}[Equivalence of Paths]
We say that $p$ is \emph{equivalent} to $p'$, abbreviated by $p \sim p'$ -- if and only if $c(p) = c(p')$ for all $c\in [\bl,\bu]$.
\end{definition}
Observe that $p \sim p'$ is not equivalent to $p\preceq p'$ and $p' \preceq p$. Instead, a criterion for $p \sim p'$ can be formulated as follows.
\begin{restatable}[]{lemma}{constantpaths}
\label{lem:constant_paths}
    Let $p, p'\in \pathst$. Then the following are equivalent:
    \begin{enumerate}
        \item $p \sim p'$.
        \item $\bestreal{p}(p') = \bestreal{p}(p)$ and $\bestreal{p'}(p) = \bestreal{p'}(p')$.
        \item $\bestreal{p}(p') \leq \bestreal{p}(p)$ and $\bestreal{p'}(p) \leq \bestreal{p'}(p')$.
    \end{enumerate}
    Moreover, equivalent paths can only differ from each other on arcs with $\ell_a = u_a$.
\end{restatable}

\begin{restatable}[]{theorem}{equivalencenonshortestweakpath}
\label{thm:equivalence_nonshortest_weak_path}
    Let $p\in \pathst$ and let $\eshpathst$ be an essential set of shortest paths. 
    The following are equivalent:
    \begin{enumerate}
        \item \label{it:1b}
            $p\in \eshpathst$.
        \item \label{it:2b}
            For all $p'\in \eshpathst$ with $p \neq p'$ holds $p' \not\preceq p$ and $p\not\sim p'$.
        \item \label{it:3b}
            For all $p' \in \eshpathst$ with $p \neq p'$ holds $\bestreal{p}(p) < \bestreal{p}(p')$.
    \end{enumerate}
    In particular, \eispp is solved by finding a scenario-covering set of weakly non-dominated and non-equivalent paths.
\end{restatable}

We conclude by highlighting a non-essential shortest path in our running example.

\begin{example}
    In \Cref{fig:running}, the path $p_5 = (a_5)$ belongs to $\eshpathst$, because $\bestreal{p_5}(p_5) = 5 < 12 = \bestreal{p_5}(p)$ for any other path $p$. The path $p_6 = (a_6)$ is a shortest path and hence in $\shpathst$, but not part of the essential set of paths, since 
    it is weakly dominated by $p_5 = (a_5)$, as 
    $\bestreal{p_6}(p_5) = 8 \leq \bestreal{p_6}(p_6) = 8$ and 
    $\bestreal{p_5}(p_5) = 5 < \bestreal{p_5}(p_6) = 12$.
\end{example}

\section{Translation into the Multi-Objective Setting}
\label{sec:mosp}

In this chapter we will show that \cispp and \eispp 
are equivalent to a specific high-dimensional multi-objective shortest path problem.

We first need to establish some basic terms. 
Let $G$ be a graph, let $w \colon \arcs(G) \to \mathbb{R}^d$ for $d \in \mathbb{N}$ be a $d$-dimensional weight function on the arcs of $G$. 
To a path $p$ in $G$ we associate the $d$-dimensional weight $w(p) \in \mathbb{R}^d$ by $w(p) \coloneqq \sum_{a\in \arcs(p)} w(a)$. 
We impose a \emph{strict partial order}, i.e., an asymmetric, transitive, irreflexive order on the set of $s$-$t$-paths $\pathst$: 
The path $p$ \emph{weakly dominates} $p'$, denoted by $p\dleq p'$, if and only if $w(p)_{j} \leq w(p')_{j}$ for all $j\in [d]$ and there exists one entry $i\in [d]$ with $w(p)_i < w(p')_i$. 
This order is called the \emph{componentwise order} in \cite{ehrgott}. 
Moreover, we also consider the \emph{strict componentwise order} on the set of $s$-$t$-paths $\pathst$: The path $p$ \emph{dominates} $p'$, in short $p \dom p'$, if and only if $w(p)_{j} < w(p')_j$ for all $j \in [d]$. 
Note that this is also a strict partial order.

\begin{definition}[{(Weakly) Efficient Paths, (Weakly) Non-dominated Points, \cite[Def 2.1]{ehrgott}}]
    Let $p$ be a $s$-$t$-path in $G$. 
    We call $p$ an \emph{efficient path} if there exists no other $s$-$t$-path $p'$ with $w(p') \dleq w(p)$. 
    If $p$ is an efficient path, then $w(p)$ is called a \emph{non-dominated point}. 
    The set of non-dominated points is called the \emph{Pareto set} and we will denote the set of efficient $s$-$t$-paths by $\effpaths$.
    Furthermore, we call $p$ a \emph{weakly efficient path} if there exists no other $s$-$t$-path in $G$ with $w(p') \dom w(p)$. 
    Similarly, $w(p)$ is called a \emph{weakly non-dominated point} if $p$ is a weakly efficient path. 
    The set of weakly efficient $s$-$t$-paths is then denoted by $\weffpaths$, and we call the set of weakly non-dominated points the \emph{weak Pareto set}.
    
    For any set of $s$-$t$-paths $P$ we define $\paretoset(P) \coloneqq \{w(p) : p\in P\}$. 
    Thus, the Pareto set will be denoted by $\paretoset(\effpaths)$, while  $\paretoset(\weffpaths)$ describes the weak Pareto set. 
\end{definition}

\begin{remark}
    We have $\effpaths \subseteq \weffpaths$ \textnormal{\cite[Eq. 2.17]{ehrgott}} and $\paretoset(\effpaths) \subseteq \paretoset(\weffpaths)$.
\end{remark}

The notion of (weakly) efficient paths leads to the following higher-dimensional analogon of the shortest path problem.

\begin{definition}[Multi-Objective Shortest Path Problem]
    For a graph $G$, a $d$-dimensional cost function on the arcs $w \colon \arcs(G) \to \mathbb{R}^d$, and a source node $s\in \nodes(G)$ the \emph{Multi-Objective Shortest Path Problem (MOSP)} is to find a set of (weakly) efficient $s$-$t$-paths for all $t \in \nodes(G)$.
\end{definition}

We are now ready to formulate the connection of the problems in \Cref{sec:ispp} with \mosp.

\begin{definition}[Induced \mosp-Instance]
    \label{def:induced-mosp}
    Let $(G, \ell, u)$ be a graph with interval costs and let $s\in \nodes(G)$ be a source node. 
    Then we call $(G, w)$ with with $2^{|\arcs(G)|}$-dimensional arc-costs 
    \begin{equation}\label{eq:weights}
            w\colon\arcs(G) \to \mathbb{R}^{\extremepoints}, \quad \text{~with~} \quad a \mapsto (c_a)_{c\in \extremepoints} 
    \end{equation}
    the \emph{\mosp instance induced by $(G, \ell, u, s)$}.
\end{definition}

\begin{remark}\label{obs:weight_to_length}
    Note that the $2^{\arcs(G)}$-dimensional weight of a path $p\in \pathst$ can then be expressed in terms of path-length $c(p)$ of the previous section, namely by 
    \begin{equation}\label{eq:weight_to_cost}
        w(p) = \sum_{a\in \arcs(p)} w(a) = (c(p))_{c\in \extremepoints}.
    \end{equation}
\end{remark}

\begin{example}
    \Cref{fig:induced_mosp} shows the construction of a \mosp instance induced by a graph with two arcs, thus leading to a \mosp instance with $d= 2^{2}$-dimensional costs.
    It illustrates the difference between the two dominance notions: $\weffpaths = \{(a_1), (a_2)\}$, while $a_1\notin \effpaths$ since $a_2 \dleq a_1$, but $a_2 \not<_D a_1$. 
\end{example}

The dominance relations in \Cref{sec:ispp} agree with the ones in this section:

\begin{figure}
    \begin{subfigure}[t]{0.47\textwidth}
        \centering
        \begin{tikzpicture}[x=1cm,y=1cm,scale=0.8]
            \tikzstyle{a} = [line width=0.8, ->, black]
            \tikzstyle{n} = [draw, circle, inner sep = 1.5, minimum width = 10, black]
            \node[n] (s) at (0,0) {$s$};
            \node[n] (t) at (2,0) {$t$};
            \draw[a] (s) to[bend left=60, looseness=1.2] node[above, midway]{$a_1 \colon [3,4]$} (t);
            \draw[a] (s) to[bend right=60, looseness=1.2] node[below,midway] {$a_2 \colon [2,3]$} (t);
        \end{tikzpicture}
        \caption{$(G, \ell, u)$ instance with labels $a \colon [\ell_a, u_a]$.}
        \label{fig:glust_inst}
    \end{subfigure}\hfill
    \begin{subfigure}[t]{0.47\textwidth}
        \centering
        
        \begin{tikzpicture}[x=1cm,y=1cm,scale=0.8]
            \tikzstyle{a} = [line width=0.8, ->, black]
            \tikzstyle{n} = [draw, circle, inner sep = 1.5, minimum width = 10, black]
            \node[n] (s) at (0,0) {$s$};
            \node[n] (t) at (2,0) {$t$};
            \draw[a] (s) to[bend left=60, looseness=1.2] node[above, midway]{$a_1 \colon (3,3,4,4)$} (t);
            \draw[a] (s) to[bend right=60, looseness=1.2] node[below,midway] {$a_2 \colon (2,3,2,3)$} (t);
        \end{tikzpicture}
        \caption{Induced $(G, w)$-instance with weight labels $w_a$.}
        \label{fig:gw_inst}
    \end{subfigure}
    \caption{Example of efficient vs. weakly efficient paths: The upper path is in $\weffpaths$ but not in $\effpaths$, and $(3,3,4,4)$ is a point in $\paretoset(\weffpaths)$ but not in $\paretoset(\effpaths)$.}
    \label{fig:induced_mosp}
\end{figure}

\begin{restatable}[]{lemma}{equivdominance}
\label{lem:equiv_dominance}
    Let $(G, \ell, u)$ be a graph with interval costs inducing the \mosp instance $(G, w)$, and let $p, p' \in \pathst$. Then
    \begin{enumerate}
        \item 
            $p\prec p' \iff p \dom p'$.
        \item 
            $p \preceq p' \iff p \dleq p'$.
    \end{enumerate} 
\end{restatable}
\begin{proof}
    By \Cref{obs:weight_to_length}, we have of $w(p)_{c} = c(p)$ for all $c\in \extremepoints$, and then obtain
    \begin{align}
        p \prec p' 
        &\iff \forall c\in [\bl,\bu]\colon \; c(p) < c(p') 
        \overset{\text{\Cref{lem:property_bestreal_dominance}}}{\iff} \forall c\in \extremepoints\colon \; c(p) < c(p') \\ 
        &\iff \forall c\in \extremepoints\colon \; w(p)_c <w(p)_c \iff p \dom p'
    \end{align} and 
    \begin{align}
        p \preceq p'    &\iff \forall c\in [\bl,\bu]\colon \; c(p) \leq c(p') \text{~and~} \exists c'\in[\bl,\bu] \colon \;c'(p) < c'(p') \\
                        &\overset{\text{\Cref{lem:bestreal_weak_dom}}}{\iff} \forall c\in \extremepoints\colon \; c(p) \leq c(p') \text{~and~} \exists c'\in\extremepoints \colon \;c'(p) < c'(p') \\ 
                        &\iff \forall c\in \extremepoints\colon \; w(p)_c \leq w(p')_c \text{~and~} \exists c'\in\extremepoints \colon \; w(p)_{c'} < w(p')_{c'} \\ 
                        &\iff p \dleq p'. \qedhere
    \end{align}
\end{proof}

Applying \cref{thm:equivalence_nonshortestpath}, \Cref{thm:equivalence_nonshortest_weak_path}  and \Cref{lem:equiv_dominance}, we obtain the following translation of \cispp and \eispp to the multi-objective setting:
\begin{theorem}\label{thm:translation-to-mosp}
Let $(G, \ell, u)$ be a graph with interval costs and $s, t \in \nodes(G)$.
\begin{enumerate}
    \item 
        We have $\shpathst = \weffpaths$. 
        In particular, \cispp is equivalent to finding the set of weakly efficient paths for the induced \mosp instance.
    \item 
        We have $\paretoset(\eshpathst) = \paretoset(\effpaths)$. 
        In particular, \eispp is equivalent to finding a minimal set $M$ of efficient paths for the induced \mosp instance with $\paretoset(M) = \paretoset(\effpaths)$.
\end{enumerate}
\end{theorem}

The following corollary states that an essential set of shortest paths in the sense of \eispp is not only minimal w.r.t.\ inclusion, but also of minimum size:

\begin{corollary}\label{cor:minimum_set_size}
    Let $\eshpathst, \eshpathst'$ be two essential shortest path sets with $\eshpathst \neq \eshpathst$. 
    Then for any $p\in \eshpathst\setminus \eshpathst'$ there exists exactly one $p'\in \eshpathst'\setminus \eshpathst$ with $p \sim p'$. 
    In particular, any set of essential shortest paths is minimum with respect to cardinality. 
\end{corollary}
\begin{proof}
    Let $p\in \eshpathst\setminus \eshpathst'$. Since $\paretoset(\eshpathst) = \paretoset(\effpaths)$, there is a $p'\in \eshpathst'$ with $w(p) = w(p')$. But then, by \Cref{obs:weight_to_length}, we have $c(p) = c(p')$ for all $c \in [\bl, \bu]$, so that $p \sim p'$. By \Cref{thm:equivalence_nonshortest_weak_path}, $p' \notin \eshpathst$, and moreover, $p'$ must be unique.
\end{proof}

\section{A Tropical Interpretation}
\label{sec:tropical}

\Cref{thm:translation-to-mosp} shows that \cispp and \eispp can be tackled in principle by \mosp algorithms. However, the dimension of the objective is exponential in the number of arcs (\Cref{def:induced-mosp}), which renders those algorithms practically infeasible. In this section, we propose to use the min-plus algebra over the tropical semiring to replace the $2^{|\arcs(G)|}$-dimensional weights by monomials in at most $|\arcs(G)|$ variables. We will then develop a label-setting Dijkstra algorithm that uses tropical polynomials as labels.

\subsection{Tropical Polynomials}

Beside the multi-criteria optimization perspective on the interval shortest path problem, there is a natural object particularly well suited to describe our desired collections of paths -- tropical polynomials.
We refer to \cite{joswig_essentials_trop,tropicalGeom_sturmfels} for more background on tropical geometry and some of its applications.
The semiring $(\mathbb{T} = \mathbb{R} \cup \{\infty\}, \oplus, \odot)$ with the operations 
\begin{equation}
    a\oplus b \coloneqq \min\{a, b\} \quad \text{~and~} \quad
    a \odot b \coloneqq a + b
\end{equation}
is called the \emph{tropical semiring}. 
Its neutral elements are $\infty$ for tropical addition $\oplus$ and $0$ for tropical multiplication $\odot$. 
Observe that it is not a ring, as there is no additive inverse.

\begin{definition}[Tropical Polynomial]
    A \emph{tropical polynomial} $\tau \in \mathbb{T}[x_1,\dots,x_m]$ is a tropical sum of finitely many \emph{tropical monomials} in the variables $x_1, \dots,x_m$ over $\mathbb{T}$, i.e.,
    \begin{equation}
        \tau(x_1,\dots, x_m) = 
        \bigoplus_{j = 1}^{k} \left(\alpha^{(j)}_0 \odot \bigodot_{i=1}^{m} x_i^{\odot \alpha^{(j)}_i}\right)
        = \min_{j \in [k]}\left\{\alpha^{(j)}_0 + \sum_{i=1}^{m} \alpha^{(j)}_i x_i \right\}, \text{~where~} \alpha_i^{(j)} \in \mathbb{T}.
    \end{equation}
\end{definition}

Consider a graph with interval costs $(G,\ell,u)$. We can model the cost of an arc $a$ as a variable $x_a$ with $\ell_a \leq x_a \leq u_a$.
A path can then be described by the tropical product of the variables associated to its arcs, more precisely:
\begin{definition}[Tropical Path Monomial]
    Let $p$ be a simple path in a graph $G$ with $\arcs(G) = \{a_1,\dots, a_m\}$. 
    We call $\tau(p) \in \mathbb{T}[x_1,\dots, x_m]$ the \emph{tropical path monomial of $p$}, where 
    \begin{equation}
        \tau(p) \coloneqq \sum_{a_i\in \arcs(p)} x_i = \bigodot_{a_i\in \arcs(p)} x_i.
    \end{equation}
\end{definition}
\begin{remark}\label{rem:path-cost-poly-eval}
Clearly, $p$ can easily be reconstructed from $\tau(p)$. Moreover, for a given scenario $c \in [\bl,\bu]$, we have $c(p) = \tau(p)(c)$, i.e., the evaluation of $\tau(p)$ at $c$.
\end{remark}

As we are interested in obtaining shortest paths, thus taking the minimum over multiple paths, tropical polynomials are an elegant way to encode this. 

\begin{definition}[Tropical $s$-$t$-Path Polynomial]
    Let $G$ be a graph and let $P$ be a collection of $s$-$t$ paths in $G$. 
    Then the \emph{tropical path polynomial} $L(P)$ is the tropical sum of the tropical path monomials of the paths in $P$:
    \begin{equation}
        L(P) \coloneqq \bigoplus_{p\in P} \tau(p) = \min \{ \tau(p) : p \in P\}.
    \end{equation}
    For the complete and essential set of shortest $s$-$t$-paths $\shpathst$ and $\eshpathst$ for $(G, \ell, u)$, we call $L(\shpathst)$ and $L(\eshpathst)$ the \emph{complete} and \emph{essential tropical $s$-$t$-path polynomial}, respectively. 
\end{definition}

Consequently, solutions to \cispp or \eispp can be described in a natural way by a single tropical path polynomial, and for a scenario $c \in [\bl,\bu]$, evaluating this polynomial at $c$ yields the cost of a shortest path w.r.t.\ $c$. In particular, $p\in \shpathst$ if and only if there is a $c \in [\bl, \bu]$ with $c(p) = L(\pathst)(c) = L(\shpathst)(c) = L(\eshpathst)(c)$.
\begin{example}
    The last column of \Cref{tab:path_table_of_running} shows the tropical path monomials for the running example of \Cref{fig:running}. 
    It also illustrates that the extension of a path by another arc can be succinctly formulated via tropical multiplication: 
    Path $p_1$ corresponds to the extension of the $s$-$v$-path $(a_1)$ (with monomial $x_{1}$) by arc $(a_2)$. 
    Using a tropical path monomial, we obtain  $\tau(p_1) = x_{1} \odot x_{2}$.
    Moreover, the tropical path polynomials of $\pathst$, $\shpathst$, and $\eshpathst$ are:
    \begin{align}
        L(\pathst) &= \bigoplus_{p\in \pathst} \tau(p) = \min\{x_{1}+x_{2},x_{1}+x_{4},x_{3}+x_{2}, x_{3}+x_{4}, x_5, x_6\},\\
        L(\shpathst) &= \bigoplus_{p\in \shpathst} \tau(p) = \min\{x_{1}+x_{4},x_{3}+x_{2}, x_{3}+x_{4}, x_5, x_6\},\\
        L(\eshpathst) &= \bigoplus_{p\in \shpathst} \tau(p) = \min\{x_{1}+x_{4},x_{3}+x_{2}, x_{3}+x_{4}, x_5\}.
    \end{align}
    One can indeed verify that the images of $L(\pathst)$, $L(\shpathst)$, and $L(\eshpathst)$ coincide over $[\bl, \bu]$.
\end{example}
We will use tropical path polynomials in an algorithm, which will make use of the dominance concepts of \Cref{sec:ispp} and \Cref{sec:mosp}. These notions extend naturally to tropical polynomials.
\begin{definition}[Dominance of Tropical Polynomials]
    Let $(G,\ell,u)$ be a graph with interval costs and $\arcs(G) = \{a_1,\dots,a_m\}$.
    Let $\tau, \tau' \in \mathbb T[x_1, \dots, x_m]$.
    
    We say that $\tau$ dominates $\tau'$, abbreviated by $\tau \prec \tau'$ if and only if for all $ c \in [\bl, \bu]$ holds $\tau(c) < \tau'(c)$. 
    Similarly, $\tau$ weakly dominates $\tau'$, $\tau \preceq \tau'$, if and only if for all $c\in [\bl,\bu]$  holds $\tau(c) \leq \tau'(c)$ and $\tau(c') < \tau'(c')$ for some $c'\in [\bl, \bu]$. 
    Finally, we add another short-hand notation, namely $\tau \precsim \tau'$, to capture that $\tau$ is either weakly dominated by or equivalent to $\tau'$.
\end{definition}

\begin{remark}
    Recall from \Cref{rem:path-cost-poly-eval} that $c(p) = \tau(p)(c)$. Hence if $p, p' \in \pathst$, then $p \prec p'$ if and only if $\tau(p) \prec \tau(p')$, and the same holds for $\preceq$ and $\precsim$.
\end{remark}

Recalling the benefits of the best scenarios from \Cref{sec:ispp}, the charm of tropical path monomials is that dominance checks can be accomplished by mere evaluations.
Identifying paths with their tropical path monomials, we have the following reformulation:

\begin{theorem}
    To solve \cispp is to find a complete tropical $s$-$t$-path polynomial, and to solve \eispp is to find an essential tropical $s$-$t$-path polynomial, for all $t\in \nodes(G)$.      
\end{theorem}

\subsection{A Tropical Dijkstra Algorithm}

We can finally proceed to present a tropical Dijkstra algorithm that solves \cispp and \eispp through tropical path polynomials.
In view of \Cref{thm:equivalence_nonshortestpath} and \Cref{thm:equivalence_nonshortest_weak_path}, the procedure for both problems differs only in the relation:
For the complete set, we use $\prec$, while for the essential set $\precsim$ is needed. 

A rough overview over our algorithm is as follows: 
We traverse the graph starting from the source node $s$, use tropical path monomials to encode any path we encounter, check if the found path is dominated by any previously found paths. 
If not, it is added by tropical addition to the corresponding tropical path polynomial, which we use to encode our set of shortest paths. 
An exact description can be found in \Cref{alg:tropdijkstra}.

\begin{algorithm}
    \caption{Tropical Dijkstra: \textcolor{darkred}{($P_{st}$, $\vartriangleleft$) } = ($\shpathst$, $\prec$) and \textcolor{darkred}{($P_{st}$, $\vartriangleleft$) } = ($\eshpathst$, $\precsim$) for the complete and essential tropical $s$-$t$-path polynomials, respectively.} 
    \label{alg:tropdijkstra}
    \begin{algorithmic}[1] 
        \Require graph with interval costs $(G, \ell, u)$, source node $s\in \nodes(G)$
        \Ensure tropical $s$-$t$-path polynomials $L_t = L(\textcolor{darkred}{P_{st}})$  for $t\in \nodes(G)$
        \State $L_v = \infty$ for all $v\in \nodes(G)$ \label{alg_it:initlabel}
        \State $L_s = 0$\label{alg_it:initstartlabel}
        \State $Q \gets $ init\_prio\_queue({$0  : (s, 0)$})\label{algit:init_queue}
        \While{$Q \neq \emptyset$}
            \State $u, \tau_u \gets$ $Q.$extract\_min() \label{alg_line:queue_extract}
            \For {$(u,v)\in \delta^+(u)$}
                \State $\tau_v \gets \tau_u \odot x_{u,v}$ \label{alg_it:path_extension}
                \If{$L_v \,{\textcolor{darkred}{\ntriangleleft}}\, \tau_v$} \label{alg_line:prec}
                    \State $prio \gets \tau_v(\bestreal{\tau_v})$ \label{alg_line:prio}
                    \State $Q$.insert($prio: (v, \tau_v)$) \label{alg_line:queue_insert}
                    \State remove all monomials $\tau_v'$ from $L_v$ with $\tau_v \,{\textcolor{darkred}{\vartriangleleft}}\, \tau_v'$ \label{alg_line:clean}
                    \State $L_v \gets L_v \oplus \tau_v$ \label{alg_line:insert}
                \EndIf \label{alg_line:endif}
            \EndFor
        \EndWhile
        \State \Return $L_t$ for $t\in \nodes(G)$
    \end{algorithmic}
\end{algorithm}
Let us explain the algorithm in more detail.
The input is a graph with interval costs $(G,\ell,u)$, as well as a source node $s \in \nodes(G)$. 
In line \eqref{alg_it:initlabel} we initialise a tropical path polynomial $L_v$ for each node $v\in \nodes(G)$ with $\infty$, essentially marking all nodes as undiscovered, while in the next line, we set $L_s$ to the constant $0$, which can be interpreted as any $s$-$s$-path having cost 0. 
We use a minimum-priority queue $Q$, which stores pairs $(v,\tau_v)$, where $v\in \nodes(G)$ and $\tau_v$ is a tropical $s$-$v$-path monomial. 
Each such tuple is inserted into the queue by the value of $\tau_v$ at $\bestreal{\tau_v}$, where the latter is the best scenario for the path that corresponds to $\tau_v$.
As $Q$ is a minimum-priority queue, this means that we propagate paths with low best scenario costs first.
This queue is initialised, and the source node with the constant monomial $\tau_s = 0$ is inserted in line \ref{algit:init_queue}.
While it is possible, we extract the node-monomial-pair $(u, \tau_u)$ with minimum priority from the queue and explore each incident edge $(u,v)$. 
We consider the unique tropical path monomial corresponding to the path obtained by extending the path of $\tau_u$ by edge $(u,v)$ in line \ref{alg_it:path_extension}. 
In line \ref{alg_line:prec} we perform the dominance check, thus verifying if there has been found a dominating $s$-$v$-path previously (for the complete case), or which weakly dominates or is equivalent (for the essential case). 
If no such path is present, we insert the newly found path in line~\ref{alg_line:queue_insert}. 
Before we add the new tropical path monomial into the tropical $s$-$v$-path polynomial in line~\ref{alg_line:insert}, we remove in line~\ref{alg_line:clean} all monomials which were discovered and stored earlier in $L_v$, and which are now (weakly) dominated.
Note that in the essential case, equivalence cannot occur at this point, as then $L_v \precsim \tau_v$ would have held in line \ref{alg_line:prec}. 
At the end, we return the tropical $s$-$v$-path polynomials found by the procedure. 
\begin{example}
    \label{ex:tropical-dijkstra}
    To see an example of how the tropical Dijkstra works, we refer to \Cref{tab:trop_dijkstra}, where one can find the values of the key variables and objects at each iteration after line \ref{alg_line:endif} of our running example. 
    Note that the order of the arcs in $\delta^+(u)$ was processed top-to-bottom as depicted in \Cref{fig:running}. 
    Observe that in iteration 3, $L_t$ contains $x_6$, which is discarded in iteration 4, when computing the essential tropical path polynomials, while it remains in $L_t$ for the complete set: 
    In this iteration, $x_5$ is considered and compared to $L_t = x_6$.
    In both cases, the if-condition in line \ref{alg_line:prec} is true, as we have $L_t  \not \prec x_5 $ and $L_t \not \precsim x_5$. 
    However, in line \ref{alg_line:clean}, $x_6$ is removed from $L_t$ when using $\precsim$, since $x_5 \precsim x_6$, yet remains in $L_t$ for $\prec$, since $x_5 \not \prec x_6$. 
    Furthermore, from iteration 6 to 7, none of the polynomials change: 
    In this case, we consider the tropical path monomial $\tau_v  = x_1 \odot x_2$ -- which we know from the first example not to be a shortest path. 
    Indeed, both $L_t \prec \tau_v$ and $L_t \precsim \tau_v$ hold. 
    There are four more iterations of the algorithm which are not depicted in the table anymore, as at the end of iteration 8 all entries in the queue are of the form $(t, \cdot)$ and node $t$ has no more outgoing edges, such that the tropical path polynomials $L_t$ and $L_v$ remain unchanged. 
\end{example}

\begin{table}
    
\setlength\tabcolsep{1.8 pt}

\begin{tabular}{|c|rrrr|rr|rr|}
\hline
$\#$ it  &\mc{1}{c|}{$1$}              &\mc{1}{c|}{$2$}           &\mc{1}{c|}{$3$}          &\mc{1}{c|}{$4$}    &\mc{1}{c|}{$5$}                &\mc{1}{c|}{$6$}    &\mc{1}{c|}{$7$}                &\mc{1}{c|}{$8$}        \\ \hline
$u$      &\mc{4}{c|}{$s$}              &\mc{2}{c|}{$v$}           &\mc{2}{c|}{$v$}                                                                                                                                          \\ \hline
$\tau_u$ &\mc{4}{c|}{$0$}              &\mc{2}{c|}{$x_3$}         &\mc{2}{c|}{$x_1$}                                                                                                                                        \\ \hline
$(u,v)$  &\mc{1}{r|}{$a_6$}            &\mc{1}{r|}{$a_1$}         &\mc{1}{r|}{$a_3$}        &$a_5$              &\mc{1}{r|}{$a_2$}              &$a_4$              &\mc{1}{r|}{$a_2$}              &$a_4$                  \\ \hline
$\tau_v$ &\mc{1}{r|}{$x_6$}            &\mc{1}{r|}{$x_1$}         &\mc{1}{r|}{$ x_3$}        &$x_5$              &\mc{1}{r|}{$x_3\odot x_2$}          &$x_3\odot x_4$          &\mc{1}{r|}{$x_1\odot x_2$}          &$x_1\odot x_4$              \\ \hline
$L_v$    &\mc{1}{r|}{$\infty$} &\mc{1}{r|}{$x_1$} &\mc{1}{r|}{$x_1$} &$x_1$       &\mc{1}{r|}{$x_1$}       &$x_1$       &\mc{1}{r|}{$x_1$}       &$x_1$           \\
         &\mc{1}{r|}{}                 &\mc{1}{r|}{}              &\mc{1}{r|}{$\oplus \; x_3$}      &$\oplus \; x_3$            &\mc{1}{r|}{$\oplus \; x_3$}            &$\oplus \; x_3$            &\mc{1}{r|}{$\oplus \; x_3$}            &$\oplus \; x_3$                \\ \hline
$L_t$    &\mc{1}{r|}{$x_6$}    &\mc{1}{r|}{$x_6$} &\mc{1}{r|}{$x_6$}&$\dr{x_6}$  &\mc{1}{r|}{$\dr{x_6}$}  &$\dr{x_6}$  &\mc{1}{r|}{$\dr{x_6}$}  &$\dr{x_6}$      \\
         &\mc{1}{r|}{}                 &\mc{1}{r|}{}              &\mc{1}{r|}{}             &$\oplus \; x_5$            &\mc{1}{r|}{$\oplus \; x_5$}             &$\oplus \; x_5$             &\mc{1}{r|}{$\oplus \; x_5$}             &$\oplus \; x_5$                 \\
         &\mc{1}{r|}{}                 &\mc{1}{r|}{}              &\mc{1}{r|}{}             &                   &\mc{1}{r|}{$\oplus \; x_3\odot x_2$}        &$x_3\odot x_2$         &\mc{1}{r|}{$\oplus \; x_3\odot x_2$}         &$\oplus \; x_3\odot x_2$             \\
         &\mc{1}{r|}{}                 &\mc{1}{r|}{}              &\mc{1}{r|}{}             &                   &\mc{1}{r|}{}                   &$\oplus \; x_3\odot x_4$        &\mc{1}{r|}{$\oplus \; x_3\odot x_4$}        &$\oplus \; x_3\odot x_4$             \\
         &\mc{1}{r|}{}                 &\mc{1}{r|}{}              &\mc{1}{r|}{}             &                   &\mc{1}{r|}{}                   &                   &\mc{1}{r|}{}                   &$\oplus \; x_1\odot x_4$            \\ \hline
$Q$      &\mc{1}{r|}{$8{:}(t,x_6)$}    &\mc{1}{r|}{$5{:}(v,x_1)$} &\mc{1}{r|}{$2{:}(v,x_3)$}&$2{:}(v,x_3)$      &\mc{1}{r|}{$5{:}(v,x_1)$}      &$5{:}(v,x_1)$      &\mc{1}{r|}{$5{:}(t,x_5)$}      &$5{:}(t,x_5)$          \\
         &\mc{1}{r|}{}                 &\mc{1}{r|}{$8{:}(t,x_6)$} &\mc{1}{r|}{$5{:}(v,x_1)$}&$5{:}(v,x_1)$      &\mc{1}{r|}{$5{:}(t,x_5)$}      &$5{:}(t,x_5)$      &\mc{1}{r|}{$6{:}(t,x_3\odot x_2)$}  &$6{:}(t,x_3\odot x_2)$      \\
         &\mc{1}{r|}{}                 &\mc{1}{r|}{}              &\mc{1}{r|}{$8{:}(t,x_6)$}&$5{:}(t,x_5)$      &\mc{1}{r|}{$6{:}(t,x_3\odot x_2)$}  &$6{:}(t,x_3\odot x_2)$  &\mc{1}{r|}{$7{:}(t,x_3\odot x_4)$}  &$7{:}(t,x_3\odot x_4)$      \\
         &\mc{1}{r|}{}                 &\mc{1}{r|}{}              &\mc{1}{r|}{}             &$8{:}(t,x_6)$      &\mc{1}{r|}{$8{:}(t,x_6)$}      &$7{:}(t,x_3\odot x_4)$  &\mc{1}{r|}{$8{:}(t,x_6)$}      &$8{:}(t,x_6)$          \\
         &\mc{1}{r|}{}                 &\mc{1}{r|}{}              &\mc{1}{r|}{}             &                   &\mc{1}{r|}{}                   &$8{:}(t,x_6)$      &\mc{1}{r|}{}                   &                       \\ \hline
\end{tabular}

    \caption{Example run of Tropical Dijkstra (\Cref{alg:tropdijkstra}) for the instance of \Cref{fig:running}: State of key variables in line \ref{alg_line:endif}, red tropical monomials are only in the complete version, black values are for both the complete and essential polynomials.}
    \label{tab:trop_dijkstra}
\end{table}

Our algorithm resembles Martins' algorithm \cite{martins_multicriteria_1984} and the MDA \cite{maristany_de_las_casas_improved_2021}, and its correctness can be proved along the same lines. \Cref{alg:tropdijkstra} can be considered as label-setting in the sense of \cite{maristany_de_las_casas_new_2024}: Although a node might be inserted into the queue multiple times, each path is discovered at most once, and can thus only appear in $Q$ at most once. A monomial may be inserted to $L_v$ and deleted later, however, the following provides a criterion to ensure permanence. In particular, once a monomial is extracted from the queue, it cannot be removed from the label anymore, so that only permanently labelled paths are propagated.

\begin{restatable}[]{lemma}{permanentlabel}
\label{lem:permanent_label}
    Let $\tau$ be a tropical path polynomial in $L_v$. Then $\tau$ cannot be removed from $L_v$ in a later iteration once the minimum priority-value of the queue $Q$ exceeds $\tau(\bestreal{\tau})$.
\end{restatable}
\begin{proof}
    This follows directly from two considerations, namely -- as in standard Dijkstra and due to the non-negative costs on the arcs -- the minimum value of the queue is monotonically increasing, meaning that once a tropical path monomial $\tau'$ with best scenario costs $\tau'(\bestreal{\tau'})$ is extracted from the queue, it is not possible to find a $\tau$ at a later iteration with lower best scenario costs.
    The second important observation is that a removal of a monomial $\tau'_v$ from $L_v$ (cf.~line \ref{alg_line:clean}) can only occur if  $\tau_v(\bestreal{\tau_v'}) \vartriangleleft \tau'_v(\bestreal{\tau_v'})$ for $\vartriangleleft\, \in \{\prec, \preceq\}$ by \Cref{lem:property_bestreal_dominance} and \Cref{lem:bestreal_weak_dom}. 
    Suppose $\tau_v$ is constructed from $\tau_u \odot x_{u,v}$ and $q_{\min}$ is the priority value by which $\tau_u$ got extracted from the queue, i.e., $\tau_u(\bestreal{\tau_u}) = q_{\min}$. 
    This then implies that $\tau_v(\bestreal{\tau_v}) \geq q_{\min}$.
    Together with $\tau_v(\bestreal{\tau_v}) \leq \tau_v(\bestreal{\tau_v'})$, we must have $\tau_v'(\bestreal{\tau_v'}) \geq q_{\min}$.

    We conclude that $\tau_v'$ is permanently in $L_v$ (i.e., cannot be dominated by some future path monomial) once the minimum value of the queue $Q$ is larger than to $\tau_v'(\bestreal{\tau_v'})$.
\end{proof}

\subsection{Implementation Details and Practical Enrichments}
There are some minor adjustments that one could think of in order to improve the algorithm. 
\begin{itemize}
    \item \textbf{Fixed Arcs}: 
        As we will see in \Cref{subsec:instances}, in our instances there are many \emph{fixed} arcs, i.e., arcs with $\ell_a = u_a$.
        On those there is only a single scenario value to choose from. 
        Moreover, our instances are structured in a way such that for a given $s$-$t$-path $p$, the set of \emph{variable} arcs, i.e., those with $\ell_a < u_a$, are in fact sufficient to determine $p$ uniquely. 
        We therefore implemented a modification of the tropical path monomials, where we extend the definition to also allow $\alpha_0 \geq 0$ and only store variables corresponding to variable arcs, while adding the fixed cost value to the $\alpha_0$ term. 
        In a post-processing step, the entire paths can then be reconstructed easily, by determining the missing sub-paths of appropriate costs. 
    
    \item \textbf{Cutoff Condition}:
        As we know, a path can only be in the set of shortest paths if it also a shortest path with respect to its best scenario. 
        Moreover, any $s$-$t$-path $p$ that is a shortest path must have costs $\bestreal{p}(p) \leq \Delta_t$, where $\Delta_t$ denotes the cost of a shortest $s$-$t$-path in $G$ w.r.t.~costs $u$.  
        We therefore propose to perform a single round of standard Dijkstra on $G$ with respect to scenario costs $u$ to obtain $\Delta_t$ for all $t\in \nodes(G)$. 
        One can then discard a path $p$ if its best scenario costs surpass $\Delta_t$, before performing the dominance check. 
        Line \ref{alg_line:prec} could then be adjusted to
        \[ \textbf{if } \tau_v(\bestreal{\tau_v}) \leq \Delta_v \textbf{ and } L_v \textcolor{darkred}{\ntriangleleft} \tau_v \textbf{ do}.\]
        If we have a guarantee that each $s$-$v$-path contains at least one variable arc, we can even use $\tau_v(\bestreal{\tau_v}) < \Delta_v$.
    
    \item \textbf{Reduced Dominance Check}:
        We tried two implementations of the tropical path polynomials: One in which we store the tropical path monomials as a vector by order of discovery, and one in which we store them in order of their best scenario costs. Using \Cref{lem:permanent_label} allows us then to compare new paths in both the dominance check (line \ref{alg_line:prec}) and the correction step (line \ref{alg_line:clean}) only with those not yet permanently labelled. 
        However, for our instances it did not pay off to maintain the sorted structure, and we used the unsorted vector implementation for our tests. 
\end{itemize}

\section{Passenger Paths in Public Transport}
\label{sec:application}
Recall that the motivation behind the examination of the interval shortest path problems was a very practical one: Is it true that only a small fraction of all paths will be shortest passenger routes, regardless of the specific choice of timetable? 
While the algorithm works in a general setting, we developed it for this specific application, meaning that we are interested in computing realistic passenger paths in a public transport network. As we already saw from the extreme case of \Cref{lem:exponential_shortest_path} the total set of shortest paths can be exponential in size, thus we do not think that experiments on generalised graphs would be particularly meaningful. Instead, we will focus our experiments only on relevant instances for our setting, namely on \emph{condensed passenger event-activity networks}, and make small adaptations to our algorithm that are appropriate for this application.

\subsection{Condensed Passenger Event-Activity Networks}
\label{subsec:condensed_pass_eans}
Event-Activity Networks are the standard approach for modelling timetables; they are directed graphs, whose nodes represent \emph{events}, such as the arrival or the departure of a vehicle at a station, while the arcs correspond to \emph{activities}, such as driving, dwelling, and transferring, which model the transition from one event to the other. For example, the departure node of one line is linked by a driving arc to the arrival node of the same line at the following station. Depending on the modelling purpose, there are multiple different incarnations of event-activity networks, modelling different aspects with different depth. Our focus is on the passenger perspective, such that our event-activity networks will consist only of passenger-related activities. Our construction is derived from a transportation network, and follows the standard approach, similar to, e.g., \cite{schiewe_periodic_2020, schmidt_integrating_2014, masing_forward}. To formally construct it, we first need to specify our input from which we derive our relevant graph with interval costs. 

We consider an undirected graph, the \emph{augmented transportation network} $\ptn$, whose nodes can be separated into a set of \emph{cells} $\Cells$ and a set of \emph{stations} $\Stations$. An edge is either an \emph{access edge} linking a cell with a station, or a \emph{vehicle edge} between two stations.
Cells represent a specific area of the city, e.g., a residential block or a section of a commercial district. 
The access edges $(c,s) \in \Edges(\ptn) \cap (\Cells \times \Stations)$ indicate the possibility of accessing a station from a cell and vice versa within a specified access time $t_{\mathrm{acc}}^{(c,s)}$.
Furthermore, we are given a \emph{line plan} $\mathcal L$, by which we mean a set of simple paths on the subgraph given by stations and vehicle edges.
Each line $l\in \mathcal L$ is associated with a frequency $f_l$, indicating how often a line is operated per period time $T\in \mathbb{N}$. 
We associate a minimum and maximum \emph{driving time} with each edge $(S,S') \in \Edges(l)$ of a line $l\in \mathcal L$, denoted as $t_{\min}^{(S,S')}$ and $t_{\max}^{(S,S')}$, as well as minimum and maximum \emph{dwelling times} $t_{\min}^S$ and $t_{\max}^S$ at each internal station $S\in \nodes(l)$ along its path. 
Finally, we are given a \emph{minimum transfer time} $t_{trans}^S$ for each station $s\in \Stations$. 

\begin{example}
    \Cref{fig:transportation_network} depicts a small example of an augmented transportation network.
    Here, it is possible to access both station $C$ and $D$ from the right cell, while the left cell has only access to $A$. 
    The green line takes $5$ minutes to drive between $A$ and $D$, and runs every $30$ minutes, while the purple line runs every $10$ minutes, corresponding to a frequency of $2$ and $6$ respectively for a period time $T = 60$. 
    \end{example}
\begin{figure}
\centering
\begin{subfigure}[t]{\textwidth}
    \centering
    \begin{tikzpicture}[xscale = 2, yscale = 2, transform shape]
    \tikzstyle{l1} = [tabpurple]
    \tikzstyle{l2} = [taborange]
    \tikzstyle{l3} = [tabgreen]
    \tikzstyle{ll1} = [ circle, l1, fill = white, font= \tiny, inner sep = 0]
    \tikzstyle{ll2} = [ circle, l2, fill = white ,font= \tiny, inner sep = 0]
    \tikzstyle{ll3} = [ circle, l3, fill = white, font= \tiny, inner sep = 0]
    \tikzstyle{lt} = [ circle, gray, fill = white, font= \tiny, inner sep = 1]
    
    \tikzstyle{station} = [draw, gray, fill = tabyellow, opacity = .5, rectangle, thick]
    \tikzstyle{line1} = [line width = 4, l1]
    \tikzstyle{line2} = [line width = 4, l2]
    \tikzstyle{line3} = [line width = 4, l3]
    
    \tikzstyle{cell} = [draw, gray, fill = gray!50]
    \tikzstyle{access} = [line width =  1.5, dotted, gray]

    \node[station] (A) at (0,0) {\scriptsize  A};
    \node[station] (B) at (2,0) {\scriptsize B};
    \node[station] (C) at (4,0) {\scriptsize C};
    \node[station] (D) at (2,2) {\scriptsize D};
    \node[cell] (c1) at (0,1) {};
    \node[cell] (c2) at (3,2) {};
    
    \draw[line1] (A) to node[ll1, midway] {$3\mid T/f{=}10$} (B);
    \draw[line1] (B) to node[ll1, midway] {$4\mid T/f{=}10$} (C);
    \draw[line2] (B) to node[ll2, sloped] {$4 \mid T/f{=}5$} (D);
    \draw[line3] (A) to node[ll3, sloped] {$5 \mid T/f{=}30$} (D);
    
    \draw[access] (c1) to node[lt, midway] {3} (A);
    \draw[access] (c2) to node[lt, midway] {4} (D); 
    \draw[access] (c2) to node[lt, midway] {9}(C);
\end{tikzpicture}
    \caption{Augmented transportation network. 
    Access edges are labelled with access times, vehicle edges with driving times and service intervals.}
    \label{fig:transportation_network}
\end{subfigure}
\hfill
\begin{subfigure}[t]{\textwidth}
    \centering
    \scalebox{0.9}{\begin{tikzpicture}[x=12mm, y=14mm]
	
	\tikzstyle{cell} = [draw, rectangle, gray, thick, fill = gray!50, minimum width = 17, minimum height = 17 ]
	\tikzstyle{a} = [line width=1, -stealth]
	\tikzstyle{transfer} = [a, dotted]
	\tikzstyle{exit} = [a, gray]
	\tikzstyle{enter} = [a, dotted]
	\tikzstyle{l1} = [tabpurple]
	\tikzstyle{l2} = [taborange]
	\tikzstyle{l3} = [tabgreen]
	\tikzstyle{border} = [a, line width=0.8, black!60, dotted]
	\tikzstyle{train} = [draw, circle, inner sep=1, minimum width=8 , scale =1, line width = 1.5]
	\tikzstyle{t1} = [train, l1]
	\tikzstyle{t2} = [ train, l2]
	\tikzstyle{t3} = [ train, l3]
	\tikzstyle{path} = [line width=4, zibLightblue]
	\tikzstyle{dep} = [train, black, fill=gray!50]
	\tikzstyle{arr} = [train, black,fill = white]
	\tikzstyle{dep1} = [t1, fill=gray!50]
	\tikzstyle{arr1} = [t1, fill=white]
	\tikzstyle{dep2} = [t2, fill=gray!50]
	\tikzstyle{arr2} = [t2, fill=white]
	\tikzstyle{dep3} = [t3, fill=gray!50]
	\tikzstyle{arr3} = [t3, fill=white]

	\tikzstyle{rop} = [opacity = 0.25]
	
	\tikzstyle{station} = [tabyellow, opacity = .8, thick  ]

    \tikzstyle{midlabel} = [ inner sep = 0, font = \tiny, sloped]
    \tikzstyle{ml1} = [midlabel, l1, fill = white]
    \tikzstyle{ml2} = [midlabel, l2, fill = white]
    \tikzstyle{ml3} = [midlabel, l3, fill = white]
    \tikzstyle{mlt} = [midlabel, gray, fill = white, sloped]
	\draw[station] (3.7,-.3) rectangle (7.3,2.3);
	\draw[station] (-.3,-.3) rectangle (2.3,2.3);
	\draw[station] (3.7,2.7) rectangle (6.3,5.3);
	\draw[station] (9.7,-.3) rectangle (10.3,1.3);

	\node[dep1] (ADH) at (0,0)  {};
	\node[arr1] (BAH) at (4,0)  {};
	\node[dep1] (BDH) at (7,0) {};
	\node[arr1] (CAH) at (10, 0) {};
	
	\node[dep1] (CDR) at (10,1)  {};
	\node[arr1] (BAR) at (7,1)  {};
	\node[dep1] (BDR) at (4,1) {};
	\node[arr1] (AAR) at (0, 1) {};
	
	\node[dep2] (DDH2) at (5, 3) {};
	\node[arr2] (BAH2) at (5, 2) {};
	
	\node[dep2] (BDR2) at (6, 2) {};
	\node[arr2] (DAR2) at (6, 3) {};

	\node[dep3] (ADH3) at (1,2) {};
	\node[arr3] (DAH3) at (4,5) {};
	
	\node[dep3] (DDR3) at (4,4) {};
	\node[arr3] (AAR3) at (2,2) {};

	\draw[a, l1] (ADH) to node[ml1] {\small [3,3]} (BAH); 
	\draw[a, l1] (BAH) to node[ml1] {\small [1,1]} (BDH); 
	\draw[a, l1] (BDH) to node[ml1] {\small [4,4]} (CAH);
	
	\draw[a, l1] (CDR) to node[ml1] {\small [4,4]} (BAR); 
	\draw[a, l1] (BAR) to node[ml1] {\small [1,1]} (BDR);
	\draw[a, l1] (BDR) to node[ml1] {\small [3,3]} (AAR); 
	
	\draw[a, l2] (DDH2) to node[ml2, sloped] {\small [4,4]} (BAH2); 
	\draw[a, l2] (BDR2) to node[ml2, sloped] {\small [4,4]} (DAR2);
	
	
	\draw[a, l3] (ADH3) to node[ml3, sloped] {\small [5,5]} (DAH3); 
	\draw[a, l3] (DDR3) to node[ml3, sloped] {\small [5,5]} (AAR3);
	
	\draw[transfer, l1] (BAH2) to node[ml1] {\small [2,11]} (BDR); 
	\draw[transfer, l2] (BAH) to node[ml2, left] {\small [2,6]} (BDR2);
	\draw[transfer, l2] (BAR) to  node[ml2] {\small [2,6]} (BDR2);
	\draw[transfer, l1] (BAH2) to node[ml1, right] {\small [2,11]} (BDH);
	\draw[transfer, l1] (AAR3) to node[ml1] {\small [2,11]} (ADH); 
	
	\draw[transfer, l3] (AAR) to node[ml3] {\small [2,31]} (ADH3);
	\draw[transfer, l3] (DAR2) to node[ml3] {\small [2,31]} (DDR3); 
	\draw[transfer, l2] (DAH3) to node[ml2, near start] {\small [2,6]} (DDH2);

	\node[cell] (s1) at (-1,2) {\small $s_1$};
	\node[cell] (t1) at (0.5,3) {\small $t_1$};
    \node[cell] (s2) at (8,5) {\small $s_2$};
	\node[cell] (t2) at (8,4) {\small $t_2$};	
	
	\draw[enter, l1] (s1) to node[ml1] {\small [3,12]} (ADH); 
	\draw[enter, l3] (s1) to node[ml3, left] {\small [3,32]} (ADH3); 
	\draw[exit] (AAR) to node[mlt,left] {\small [3,3]} (t1); 
	\draw[exit] (AAR3) to node[mlt] {\small [3,3]} (t1); 
	\draw[exit] (DAR2) to node[mlt] {\small [4,4]}  (t2); 
	\draw[exit] (DAH3) to node[mlt, left, near start] {\small [4,4]} (t2); 
	\draw[enter, l2] (s2) to node[ml2,right] {\small [4,8]} (DDH2); 
	\draw[enter, l3] (s2) to node [ml3, near start] {\small [4,33]}(DDR3);
	\draw[exit] (CAH) to node[mlt] {\small [9,9]} (t2); 
	\draw[enter,l1] (s2) to node[ml1] {\small [9,19]} (CDR); 

\end{tikzpicture}}
    
    \caption{Derived event-activity network from \Cref{fig:transportation_network} with a minimum transfer time of $2$.}
    \label{fig:ean}
\end{subfigure}
\caption{Augmented transportation network and derived event-activity network. Stations are yellow, cells are grey.}
\end{figure}

As mentioned, the construction of the event-activity network follows the standard approach. We include it for completeness:

\begin{definition}[Condensed Passenger Event-Activity Network]
\label{def:condensed_pass_networks}
Given a $\ptn$, a period time $T$, a line plan $\mathcal L$, frequency assignment $f_l$ for $l\in \mathcal L$, as well as specified time estimates as above, we construct the \emph{condensed passenger event-activity network} $(G, \ell, u)$ as a graph with interval costs as follows:
\begin{itemize}
    \item For each line $l\in \mathcal L$: 
    \begin{itemize}
        \item add the \emph{driving arcs} $a= ((\dep, l, S, dir), (\arr, l, S', dir))$, for each direction $dir \in \{\vartriangleleft, \vartriangleright\}$ of each edge $(S,S') \in \Edges(l)$, with bounds $\ell_{a} \coloneqq t_{\min}^{(S,S')}$ and $u_a \coloneqq t_{\max}^{(S,S')}$.
        \item add the \emph{dwelling arcs} $a= ((\arr, l, S, dir),(\dep, l, S, dir))$, for each direction $dir \in \{\vartriangleleft, \vartriangleright\}$ for each interior node $S\in \nodes(l)$, with bounds $\ell_{a} \coloneqq t_{\min}^S$ and $u_{a} \coloneqq t_{\max}^S$.
    \end{itemize}
    We call nodes of the form $(arr, \cdot, \cdot, \cdot)$ arrival events, and $(dep, \cdot, \cdot, \cdot)$ departure events, and we denote by $Arr(S)$ and $Dep(S)$ the set of arrival and departure nodes at station $S\in \Stations$.
    \item For each $a\in Arr(S) \times Dep(S)$ with $l\neq l'$, add a \emph{transfer arc} $a = ((\arr, l, S, dir), (\dep, l', S, dir'))$, with bounds $\ell_{a} \coloneqq t_{trans}^S$ and $u_{a} \coloneqq t_{trans}^S + T/f_{l'}-1$.
    \item For each access-edge $(c, S)\in \Edges(\ptn) \cap (\Cells \times \Stations)$:
    \begin{itemize}
        \item add an \emph{access arc} $(s_c, x)$ for each departure event $x\in Dep(S)$, with bounds $\ell_{(s_c,x)} \coloneqq t_{acc}^{(c,S)}$ and $u_{(s_c,x)} \coloneqq t_{acc}^{(c,S)} + T/f_l - 1$, where $f_l$ is the frequency of the line of departure $x = (dep, l, S, dir)$.
        \item add an \emph{exit arc} $(x, t_c)$ for each arrival event $x\in Arr(S)$, with bounds both set to the access time, i.e.,  $\ell_{(x,t_c)} \coloneqq u_{(x, t_c)} \coloneqq t_{acc}^{(c, S)}$.
    \end{itemize}
    We call $s_c$ a  \emph{source cell}, and $t_c$ a \emph{target cell}.
    
\end{itemize}
\end{definition}

While the reasoning behind the bounds on the arcs should be intuitive for driving and dwelling arcs, the interval values for the transfer, access and exit arcs deserve a little discussion:  
Suppose $t$ is the access time or minimum transfer time. 
Then a passenger arriving at the location of the next departure can either be lucky and the vehicle departs right away, or it might at most take $T/f - 1$ more minutes for the next vehicle to depart, where $f$ is the frequency of the line of the corresponding departure event.   
Thus, we assign the bound $[t, t + T/f - 1]$ for the transfer and access arcs.
The exit arcs $a$ have $\ell_a = u_a$ corresponding to the access time of a station to its cell, as a passenger departing a vehicle can directly walk to their cell-location, without further delays.

One could argue that, similarly to exit activities, access activities should be modelled with fixed times as well, as in practice, a passenger would look up the timetable beforehand and go to the station ``just in time'' and would then not have to wait. 
In this case, however, this is not a good choice, as slower connections at a higher frequency might be more desirable to a passenger than a low frequency line with shorter travel times, which we also show in the next example.

\begin{example}
\label{ex:frequency_transfer}
    An example of the derived event-activity network from the augmented transportation network from \Cref{fig:transportation_network} can be found in \Cref{fig:ean}.
     The access time from the left cell to station $A$ is $3$ minutes. If a passenger arrives from this cell, they must wait less than $10$ minutes for the purple line to depart, as it runs with a frequency of $6$. Thus, the bounds of on the purple dotted arc incident to $s_1$ is $[3,12]$, while the bounds on its green counterpart are $[3,32]$, as the green line runs only every 30 minutes. 
    To see the advantage of modelling access with times dependent on the frequency of the departing lines, consider the following: If the access time at $s_1$ was fixed to $[3,3]$ each, then the shortest path from $s_1$ to $t_2$ would always be via the green line. This line runs only every 30 minutes however. Thus, if the green line has just left, the connection via the purple, or even the purple and orange line can be good alternatives, such that a passenger might arrive at $t_2$ earlier in comparison to waiting for the next trip of the green line. 
\end{example}

Note that -- unlike for timetabling purposes, where one would typically create one representative per frequency of each line -- we construct only one sequence of dwell and driving activities per line (and direction). 
This means that we do not differentiate between one frequency of one line and the next, the advantage of higher frequency service is only reflected in the reduced transfer times on the corresponding arcs. 
In this sense, it captures all relevant information for the passengers, but can be considered a condensed version of the network needed in the context of timetabling, in which higher frequencies are modelled by copies of the vehicle-related activities. 

\subsection{Algorithm Adjustments for Realistic Passenger Paths}
\label{subsec:realistic_adjustments}
The algorithm we proposed is formulated for the pure problem with a generic graph with interval costs. For the application of modelling passenger paths however, some small adjustments are desirable: 
\begin{itemize}
    \item 
\textbf{Simple Station Paths:}
While the pure algorithm will only return simple paths in the graph (under the assumption that there are no zero-cost cycles), the structure of our instances allows for `unrealistic' paths, in which passengers visit a station more than once, as a simple path in $G$ might project to a non-simple path in the augmented public transportation network $\ptn$. 
When constructing the path, i.e., in line \eqref{alg_it:path_extension} we verify that the station sequence of the path is simple. More precisely, we do so by storing the set of previously visited stations, and whenever we extend the path by a driving or exit arc, we check that the station $S$ of the tail of the arc has not been added before. 

\item \textbf{Output Restriction and Preprocessing:}
We are only interested in passenger paths, i.e., paths which start in a source and end in a target cell, thus we restrict our output to the set of output cells. 
This also allows us to do some preprocessing, as we can contract any nodes with in- and out-degree of one, whenever they are not a target cell, as then either both or none of the two succeeding arcs may be part of a passenger path. 
Moreover, we delete any target cells which are directly connected to a station node, which is also accessible from the input source node. This way, any passenger path uses at least one line of the public transport network. 

\item \textbf{Restriction of Number of Transfers:}
The number of transfers often plays a significant role for the quality of a passenger path. 
From a practical point of view, it might be desirable to allow only for paths with a limited number of transfers. 
In our case, our algorithm can be adapted to this purpose in an elegant way, without extending the data structure of the tropical polynomials: 
In our instances, $\ell_a = u_a$ holds for driving and dwelling activities $a$. As mentioned, we replace the corresponding variables in tropical path polynomials by absolute numbers, so that variables only remain for transfer and access activities. Consequently, for each path from a source to a target event with $k$ transfers, the number of variables (the tropical degree) in its tropical path monomial is $k+1$. 
When restricting the number of transfers to $k$, we thus simply check whether the degree of the tropical path monomial in line \eqref{alg_it:path_extension} does not exceed $k+1$. 
For more general instances, where also dwelling and driving activities are variable, one could easily include this as well, but would need to keep track of the number of transfers explicitly for each path monomial. 

\end{itemize}

\section{Experiments}
\label{sec:experiments}
\label{subsec:instances}
We perform experiments on two sets of realistic instances. 
The first set of 17 instances we construct from augmented transportation networks of the city of Wuppertal, following \Cref{def:condensed_pass_networks} with data provided by WSW Wuppertaler Stadtwerke GmbH. To diversify our test instances, we also include a second set, which we derive from the \timpasslib{} instances, a benchmark library for the Integrated Periodic Timetabling and Passenger Routing Problem \cite{schiewe_introducing_2023}. 

The instances of \timpasslib{} provide an event-activity network, background information about each of the events (e.g., line and station association) as well as a demand matrix. Note that the event-activity network of this benchmark library is not entirely in the format we need, e.g., there are  $f_l$ copies of each vehicle arc for each line with frequency $f_l$ needed for timetabling, as well as other auxiliary arcs, such as synchronisation and headway arcs. 
However, the instances can be seen as realistic representations of public transportation networks, and we can easily derive the condensed passenger event-activity networks as we need them for our application as follows: 

Each of the events in the event-activity network of \timpasslib{} has an associated station, line, direction and frequency-repetition value. 
The condensed passenger event activity network then contains all nodes of the original instance with frequency-repetition value $1$, and all driving, dwelling and transfer arcs between them. To each driving and dwelling arc, we assign $l_a = u_a$ as the lower bound of the corresponding activities of the original instance. 
We also use the original lower bound on the transfer arcs, but set the upper bound to $T/f_l + \ell_a -1$, where $a = (v,w)$ is the corresponding transfer arc and $f_l$ is the frequency of line $l$ of the associated event $w$, i.e., maximum frequency repetition of copies of the same event in the original. This agrees with our reasoning of transfer time bounds described in the construction of condensed passenger event-activity networks (cf. \Cref{subsec:condensed_pass_eans}). 
While \timpasslib{} does not contain information on cells, it provides demand between origin-destination pairs on the station-level. 
We interpret each origin and each destination station as a separate cell, and link them by access and exit arcs to each event of the same station. 
In both cases, we set the lower bounds to $0$ and the upper bounds to $T/f_l -1$ and $0$ for access and exit arcs respectively, where $f_l$ is the frequency of line $l$ of event $w$ for the access arc $(v,w)$. Note that, again, this agrees with the bounds on transfer times of our original construction. 

An overview of the instances -- $17$ Wuppertal instances $\Wup{1}-\Wup{17}$ and $26$ \timpasslib{} instances -- can be found in \Cref{tab:instance_overview}. Note that we do not include the \timpasslib{} instance \texttt{Erding-NDP-S021}, as its condensed passenger event-activity network coincides with that of \texttt{Erding-NDP-S020}. 
\begin{table}

	\setlength\tabcolsep{5 pt}
	\begin{tabular}{| l | r r r r r | r r r| r r}
 instance & $ \mu(G) $ &  $ |\nodes(G)| $  & $ | \arcs(G) | $  & $ | \arcs_{\mathrm{trf}} | $  &  $ | \arcs_{\mathrm{fix}} | $  & $ | \mathcal L| $ & $ | \mathcal S |$ & $ |\mathcal C| $ & $ T $ &  $u_{max}$ \\ \hline
$ \Wup{1} $ & $ 59 $ & $ 84 $ & $ 142 $ & $ 0 $ & $ 84 $ & $ 1 $ & $ 8 $ & $ 28 $ & $ 60 $ & $ 0 $ \\ 
$ \Wup{2} $ & $ 121 $ & $ 124 $ & $ 244 $ & $ 8 $ & $ 142 $ & $ 2 $ & $ 13 $ & $ 36 $ & $ 60 $ & $ 11 $ \\ 
$ \Wup{3} $ & $ 285 $ & $ 174 $ & $ 458 $ & $ 22 $ & $ 250 $ & $ 3 $ & $ 16 $ & $ 52 $ & $ 60 $ & $ 21 $ \\ 
$ \Wup{4} $ & $ 445 $ & $ 282 $ & $ 726 $ & $ 54 $ & $ 392 $ & $ 5 $ & $ 28 $ & $ 80 $ & $ 60 $ & $ 21 $ \\ 
$ \Wup{5} $ & $ 596 $ & $ 326 $ & $ 921 $ & $ 107 $ & $ 478 $ & $ 8 $ & $ 30 $ & $ 84 $ & $ 60 $ & $ 21 $ \\ 
$ \Wup{6} $ & $ 820 $ & $ 414 $ & $ 1233 $ & $ 183 $ & $ 620 $ & $ 10 $ & $ 39 $ & $ 102 $ & $ 60 $ & $ 21 $ \\ 
$ \Wup{7} $ & $ 1189 $ & $ 492 $ & $ 1680 $ & $ 324 $ & $ 800 $ & $ 13 $ & $ 44 $ & $ 111 $ & $ 60 $ & $ 21 $ \\ 
$ \Wup{8} $ & $ 1880 $ & $ 620 $ & $ 2499 $ & $ 627 $ & $ 1101 $ & $ 16 $ & $ 56 $ & $ 129 $ & $ 60 $ & $ 21 $ \\ 
$ \Wup{9} $ & $ 2218 $ & $ 694 $ & $ 2911 $ & $ 803 $ & $ 1243 $ & $ 18 $ & $ 62 $ & $ 140 $ & $ 60 $ & $ 21 $ \\ 
$ \Wup{10} $ & $ 2870 $ & $ 826 $ & $ 3695 $ & $ 1147 $ & $ 1511 $ & $ 24 $ & $ 74 $ & $ 152 $ & $ 60 $ & $ 61 $ \\ 
$ \Wup{11} $ & $ 3461 $ & $ 916 $ & $ 4376 $ & $ 1512 $ & $ 1699 $ & $ 28 $ & $ 82 $ & $ 163 $ & $ 60 $ & $ 61 $ \\ 
$ \Wup{12} $ & $ 4273 $ & $ 1030 $ & $ 5302 $ & $ 2000 $ & $ 1957 $ & $ 31 $ & $ 94 $ & $ 178 $ & $ 60 $ & $ 61 $ \\ 
$ \Wup{13} $ & $ 5482 $ & $ 1232 $ & $ 6713 $ & $ 2777 $ & $ 2347 $ & $ 37 $ & $ 107 $ & $ 200 $ & $ 60 $ & $ 61 $ \\ 
$ \Wup{14} $ & $ 6089 $ & $ 1406 $ & $ 7494 $ & $ 3072 $ & $ 2652 $ & $ 47 $ & $ 128 $ & $ 215 $ & $ 60 $ & $ 61 $ \\ 
$ \Wup{15} $ & $ 7143 $ & $ 1584 $ & $ 8726 $ & $ 3550 $ & $ 3098 $ & $ 57 $ & $ 138 $ & $ 225 $ & $ 60 $ & $ 61 $ \\ 
$ \Wup{16} $ & $ 8622 $ & $ 1752 $ & $ 10373 $ & $ 4561 $ & $ 3486 $ & $ 68 $ & $ 146 $ & $ 229 $ & $ 60 $ & $ 61 $ \\ 
$ \Wup{17} $ & $ 11384 $ & $ 1904 $ & $ 13287 $ & $ 6833 $ & $ 3874 $ & $ 77 $ & $ 148 $ & $ 229 $ & $ 60 $ & $ 61 $ \\ 
\hline
 $ \texttt{toy} $ &  $ 29 $ &  $ 48 $ &  $ 76 $ &  $ 16 $ &  $ 44 $ &  $ 2 $ &  $ 8 $ &  $ 8 $ &  $ 60 $ &  $ 62 $ \\ 
 $ \texttt{toy-2} $ &  $ 189 $ &  $ 80 $ &  $ 268 $ &  $ 152 $ &  $ 84 $ &  $ 6 $ &  $ 8 $ &  $ 8 $ &  $ 60 $ &  $ 62 $ \\ 
 $ \texttt{metro} $ &  $ 289 $ &  $ 358 $ &  $ 646 $ &  $ 142 $ &  $ 376 $ &  $ 4 $ &  $ 51 $ &  $ 51 $ &  $ 150 $ &  $ 59 $ \\ 
 $ \texttt{regional} $ &  $ 513 $ &  $ 290 $ &  $ 802 $ &  $ 376 $ &  $ 323 $ &  $ 8 $ &  $ 34 $ &  $ 27 $ &  $ 60 $ &  $ 60 $ \\ 
 $ \texttt{grid} $ &  $ 743 $ &  $ 266 $ &  $ 1008 $ &  $ 592 $ &  $ 308 $ &  $ 8 $ &  $ 25 $ &  $ 25 $ &  $ 60 $ &  $ 62 $ \\ 
 $ \texttt{Hamburg} $ &  $ 629 $ &  $ 642 $ &  $ 1270 $ &  $ 284 $ &  $ 740 $ &  $ 7 $ &  $ 68 $ &  $ 67 $ &  $ 10 $ &  $ 11 $ \\ 
 $ \texttt{Erding-NDP-S020} $ &  $ 1159 $ &  $ 548 $ &  $ 1706 $ &  $ 978 $ &  $ 589 $ &  $ 21 $ &  $ 51 $ &  $ 28 $ &  $ 60 $ &  $ 62 $ \\ 
 $ \texttt{Schweiz-Fernverkehr} $ &  $ 5260 $ &  $ 1520 $ &  $ 6779 $ &  $ 4401 $ &  $ 1773 $ &  $ 40 $ &  $ 140 $ &  $ 136 $ &  $ 120 $ &  $ 126 $ \\ 
 $ \texttt{R1L1} $ &  $ 5882 $ &  $ 4164 $ &  $ 10045 $ &  $ 2827 $ &  $ 5386 $ &  $ 55 $ &  $ 250 $ &  $ 250 $ &  $ 60 $ &  $ 62 $ \\ 
 $ \texttt{R1L2} $ &  $ 6044 $ &  $ 4168 $ &  $ 10211 $ &  $ 2983 $ &  $ 5394 $ &  $ 54 $ &  $ 250 $ &  $ 250 $ &  $ 60 $ &  $ 62 $ \\ 
 $ \texttt{R1L3} $ &  $ 6526 $ &  $ 4684 $ &  $ 11209 $ &  $ 2971 $ &  $ 6146 $ &  $ 65 $ &  $ 250 $ &  $ 250 $ &  $ 60 $ &  $ 62 $ \\ 
 $ \texttt{R2L1} $ &  $ 6797 $ &  $ 4716 $ &  $ 11512 $ &  $ 3332 $ &  $ 6102 $ &  $ 66 $ &  $ 280 $ &  $ 280 $ &  $ 60 $ &  $ 62 $ \\ 
 $ \texttt{R2L2} $ &  $ 7002 $ &  $ 4764 $ &  $ 11765 $ &  $ 3489 $ &  $ 6174 $ &  $ 66 $ &  $ 280 $ &  $ 280 $ &  $ 60 $ &  $ 62 $ \\ 
 $ \texttt{R1L4} $ &  $ 8029 $ &  $ 5260 $ &  $ 13288 $ &  $ 3910 $ &  $ 6998 $ &  $ 71 $ &  $ 250 $ &  $ 250 $ &  $ 60 $ &  $ 62 $ \\ 
 $ \texttt{R2L3} $ &  $ 7726 $ &  $ 5608 $ &  $ 13333 $ &  $ 3397 $ &  $ 7412 $ &  $ 80 $ &  $ 280 $ &  $ 280 $ &  $ 60 $ &  $ 62 $ \\ 
 $ \texttt{R3L1} $ &  $ 8549 $ &  $ 5108 $ &  $ 13656 $ &  $ 4770 $ &  $ 6628 $ &  $ 73 $ &  $ 296 $ &  $ 296 $ &  $ 60 $ &  $ 62 $ \\ 
 $ \texttt{R3L2} $ &  $ 8658 $ &  $ 5044 $ &  $ 13701 $ &  $ 4937 $ &  $ 6538 $ &  $ 70 $ &  $ 296 $ &  $ 296 $ &  $ 60 $ &  $ 62 $ \\ 
 $ \texttt{R4L1} $ &  $ 9625 $ &  $ 5570 $ &  $ 15194 $ &  $ 5502 $ &  $ 7226 $ &  $ 86 $ &  $ 319 $ &  $ 319 $ &  $ 60 $ &  $ 62 $ \\ 
 $ \texttt{R4L2} $ &  $ 10096 $ &  $ 5686 $ &  $ 15781 $ &  $ 5863 $ &  $ 7394 $ &  $ 89 $ &  $ 319 $ &  $ 319 $ &  $ 60 $ &  $ 62 $ \\ 
 $ \texttt{R3L3} $ &  $ 10574 $ &  $ 6316 $ &  $ 16889 $ &  $ 5631 $ &  $ 8396 $ &  $ 95 $ &  $ 296 $ &  $ 296 $ &  $ 60 $ &  $ 62 $ \\ 
 $ \texttt{long-distance} $ &  $ 14419 $ &  $ 3280 $ &  $ 17698 $ &  $ 12242 $ &  $ 4086 $ &  $ 42 $ &  $ 250 $ &  $ 240 $ &  $ 60 $ &  $ 62 $ \\ 
 $ \texttt{Stuttgart} $ &  $ 12937 $ &  $ 5478 $ &  $ 18414 $ &  $ 10572 $ &  $ 6113 $ &  $ 156 $ &  $ 560 $ &  $ 391 $ &  $ 3600 $ &  $ 1979 $ \\ 
 $ \texttt{R4L3} $ &  $ 12599 $ &  $ 7006 $ &  $ 19604 $ &  $ 7098 $ &  $ 9322 $ &  $ 115 $ &  $ 319 $ &  $ 319 $ &  $ 60 $ &  $ 62 $ \\ 
 $ \texttt{R2L4} $ &  $ 12609 $ &  $ 8220 $ &  $ 20828 $ &  $ 5740 $ &  $ 11258 $ &  $ 116 $ &  $ 280 $ &  $ 280 $ &  $ 60 $ &  $ 62 $ \\ 
 $ \texttt{R3L4} $ &  $ 15064 $ &  $ 8772 $ &  $ 23835 $ &  $ 7715 $ &  $ 12030 $ &  $ 120 $ &  $ 296 $ &  $ 296 $ &  $ 60 $ &  $ 62 $ \\ 
 $ \texttt{R4L4} $ &  $ 17116 $ &  $ 9022 $ &  $ 26137 $ &  $ 9635 $ &  $ 12310 $ &  $ 133 $ &  $ 319 $ &  $ 319 $ &  $ 60 $ &  $ 62 $ \\ 
    \end{tabular}
    
\caption{Overview of size metrics for our instances: cyclomatic number, number of nodes, arcs, transfer and fixed arcs, lines, stations, cells, period time, and maximal upper bound on transfer arcs. Instances are sorted by number of arcs.}
\label{tab:instance_overview}
\end{table}

\subsection{Tests and Results}

We have implemented \Cref{alg:tropdijkstra} with the practical adjustments described in \Cref{subsec:realistic_adjustments} in C++ and run our tests on an Intel i7-9700K CPU with 64 GB RAM with a 24h time limit.

Our tests are designed to answer multiple questions:
\begin{itemize}
    \item Does tropical Dijkstra compute the desired sets in a reasonable amount of time?
    \item Is our intuition correct that the number of shortest paths is indeed significantly smaller than the total number of paths?
    \item How large is the difference between the complete and the efficient set of shortest paths? 
    \item How does restricting transfers affect both performance and number of paths?
\end{itemize}
To answer these questions, we run two rows of tests:
Firstly, we compute both an essential and the complete set of shortest paths without any restrictions on the number of transfers in a smaller test setup. Secondly, we restrict to paths with at most three transfers for all instances.

\subsubsection{Paths with Unbounded Transfers}
We run both the complete and the essential version of the tropical Dijkstra on the instances $\Wup{1} - \Wup{12}$, as well as the first seven \timpasslib-instances for $10\%$ of the cells of each instance. To obtain a representative sample, we select the cells at random, but in the case of the Wuppertal instances we ensure that at least two cells are on the city outskirts and one cell has direct access to the central station. 

This initial, small set of tests is aimed at determining how well our algorithm works in the unrestricted case, at finding its limits, and at figuring out whether the obvious intuition is true, that only a small fraction of the paths are also shortest on realistic instances.

\begin{table}
    \centering
    
\setlength\tabcolsep{1.5 pt}
    \begin{tabular}{l | r r |r r |r r |r r r | r r}
 
& \multicolumn{2}{c|}{All paths ($\pathst$)} & \multicolumn{2}{c|}{Compl. ($\shpathst, \prec$)}  & \multicolumn{2}{c|}{Essential ($\eshpathst, \precsim$)} & \multicolumn{3}{c|}{Ratio/st} & \multicolumn{2}{c}{Time/s} \\
inst &  P/s & P/st &  P/s & P/st &  P/s & P/st & $ \frac{\shpathst}{\pathst} $ & $ \frac{\eshpathst}{\pathst} $  & $\frac{\eshpathst}{\shpathst}$ & $\prec$ & $\precsim$ \\ \hline
 $ \texttt{W1} $ &  $ 32.00 $  &  $ 1.38 $  &  $ 30.67 $  &  $ 1.33 $  &  $ 30.67 $  &  $ 1.33 $  &  $ 0.96 $  &  $ 0.96 $  &  $ 1.00 $  &  $ 0.00 \mathrm{ms} $  &  $ 0.00 \mathrm{ms} $  \\ 
 $ \texttt{W2} $ &  $ 60.75 $  &  $ 1.94 $  &  $ 54.75 $  &  $ 1.75 $  &  $ 53.00 $  &  $ 1.70 $  &  $ 0.90 $  &  $ 0.87 $  &  $ 0.97 $  &  $ 0.00 \mathrm{ms} $  &  $ 0.00 \mathrm{ms} $  \\ 
 $ \texttt{W3} $ &  $ 128.33 $  &  $ 2.89 $  &  $ 95.83 $  &  $ 2.16 $  &  $ 95.00 $  &  $ 2.14 $  &  $ 0.75 $  &  $ 0.74 $  &  $ 0.99 $  &  $ 0.00 \mathrm{ms} $  &  $ 0.00 \mathrm{ms} $  \\ 
 $ \texttt{W4} $ &  $ 205.25 $  &  $ 2.73 $  &  $ 155.12 $  &  $ 2.06 $  &  $ 152.38 $  &  $ 2.02 $  &  $ 0.76 $  &  $ 0.74 $  &  $ 0.98 $  &  $ 0.00 \mathrm{ms} $  &  $ 0.00 \mathrm{ms} $  \\ 
 $ \texttt{W5} $ &  $ 1768.89 $  &  $ 22.65 $  &  $ 577.11 $  &  $ 7.38 $  &  $ 557.56 $  &  $ 7.13 $  &  $ 0.33 $  &  $ 0.32 $  &  $ 0.97 $  &  $ 3.89 \mathrm{ms} $  &  $ 3.89 \mathrm{ms} $  \\ 
 $ \texttt{W6} $ &  $ 6911.64 $  &  $ 71.56 $  &  $ 741.64 $  &  $ 7.68 $  &  $ 682.00 $  &  $ 7.06 $  &  $ 0.11 $  &  $ 0.10 $  &  $ 0.92 $  &  $ 9.82 \mathrm{ms} $  &  $ 9.09 \mathrm{ms} $  \\ 
 $ \texttt{W7} $ & 1.06e+05 &  $ 1013.48 $  &  $ 2510.42 $  &  $ 24.01 $  &  $ 2241.17 $  &  $ 21.44 $  &  $ 0.02 $  &  $ 0.02 $  &  $ 0.89 $  &  $ 68.92 \mathrm{ms} $  &  $ 58.08 \mathrm{ms} $  \\ 
 $ \texttt{W8} $ & 2.03e+05 &  $ 1675.49 $  & 1.34e+04 &  $ 110.69 $  & 1.21e+04 &  $ 99.75 $  &  $ 0.07 $  &  $ 0.06 $  &  $ 0.90 $  &  $ 4.63 \mathrm{s} $  &  $ 4.04 \mathrm{s} $  \\ 
 $ \texttt{W9} $ & 1.39e+06 & 1.03e+04 & 2.67e+04 &  $ 198.12 $  & 2.09e+04 &  $ 155.47 $  &  $ 0.02 $  &  $ 0.02 $  &  $ 0.78 $  &  $ 14.22 \mathrm{s} $  &  $ 8.81 \mathrm{s} $  \\ 
 $ \texttt{W10} $ &   &   & 4.55e+04 &  $ 308.16 $  & 3.61e+04 &  $ 244.40 $  &   &   &  $ 0.79 $  &  $ 91.14 \mathrm{s} $  &  $ 61.82 \mathrm{s} $  \\ 
 $ \texttt{W11} $ &   &   & 6.58e+04 &  $ 416.92 $  & 5.48e+04 &  $ 347.41 $  &   &   &  $ 0.83 $  &  $ 275.84 \mathrm{s} $  &  $ 209.93 \mathrm{s} $  \\ 
 $ \texttt{W12} $ &   &   & 1.30e+05 &  $ 764.38 $  & 1.07e+05 &  $ 627.65 $  &   &   &  $ 0.82 $  &  $ 933.48 \mathrm{s} $  &  $ 637.46 \mathrm{s} $  \\ 
\hline
 $ \texttt{toy} $ &  $ 14.00 $  &  $ 2.00 $  &  $ 14.00 $  &  $ 2.00 $  &  $ 14.00 $  &  $ 2.00 $  &  $ 1.00 $  &  $ 1.00 $  &  $ 1.00 $  &  $ 0.00 \mathrm{ms} $  &  $ 0.00 \mathrm{ms} $  \\ 
 $ \texttt{toy-2} $ &  $ 684.00 $  &  $ 97.71 $  &  $ 118.00 $  &  $ 16.86 $  &  $ 118.00 $  &  $ 16.86 $  &  $ 0.17 $  &  $ 0.17 $  &  $ 1.00 $  &  $ 8.00 \mathrm{ms} $  &  $ 8.00 \mathrm{ms} $  \\ 
 $ \texttt{metro} $ &  $ 6388.00 $  &  $ 127.76 $  &  $ 3701.33 $  &  $ 74.03 $  &  $ 185.67 $  &  $ 3.71 $  &  $ 0.58 $  &  $ 0.03 $  &  $ 0.05 $  &  $ 540.67 \mathrm{ms} $  &  $ 5.50 \mathrm{ms} $  \\ 
 $ \texttt{regional} $ & 2.25e+05 &  $ 8662.21 $  & 3.94e+04 &  $ 1516.97 $  &  $ 7070.67 $  &  $ 271.95 $  &  $ 0.18 $  &  $ 0.03 $  &  $ 0.18 $  &  $ 10916.43 \mathrm{s} $  &  $ 207.73 \mathrm{s} $  \\ 
 $ \texttt{grid} $ & 8.21e+04 &  $ 3421.82 $  &  $ 835.33 $  &  $ 34.81 $  &  $ 835.33 $  &  $ 34.81 $  &  $ 0.01 $  &  $ 0.01 $  &  $ 1.00 $  &  $ 3.92 \mathrm{s} $  &  $ 3.92 \mathrm{s} $  \\ 
 $ \texttt{Hamburg} $ &  $ 5926.14 $  &  $ 89.79 $  &  $ 1133.43 $  &  $ 17.17 $  &  $ 1071.14 $  &  $ 16.23 $  &  $ 0.19 $  &  $ 0.18 $  &  $ 0.95 $  &  $ 31.71 \mathrm{ms} $  &  $ 27.14 \mathrm{ms} $  \\ 
   $ \texttt{Erding-NDP-S020}$  &  $ - $  &  $ - $  &  $ - $  &  $ - $  &  $ - $  &  $ - $  &  $ - $  &  $ - $  &  $ - $ &  $ - $   &  $ - $ \\ 
\end{tabular}

    \caption{Average values of our computations with no restriction on the number of transfers. P/s and P/st indicate the average number of paths found per source node and per OD pair, for the set of all paths, and the complete and essential set returned by tropical Dijkstra. Columns 7-9 show the ratios between the sizes of all three sets, the average being taken over all OD pairs again. Finally, the last two columns report the average run time per source node for both dominance relations. For example, \Wup{5} has almost $22.65$ paths on average per OD pair, but only $7.38$ potentially shortest paths, and $7.13$ essential paths. This means that only $33\%$ of all paths are ever a shortest path, and 97\% of the shortest paths are essential.}
    \label{tab:all_tests_average}
\end{table}

We are able to compute both the complete and essential set of shortest paths for the instances up to \Wup{12} within a reasonable amount of time, as well as for the smallest six instances of \timpasslib. The last instance of this test row, \texttt{Erding-NDP-S020}, will be discussed later. 

As expected, the time drastically increases with the instance size, where small instances finish in well under a second, while the largest Wuppertal instance  took approximately $15$~minutes on average per source node to compute the complete set.
In comparison to the complete set computations, we observe that indeed, the runtimes on average are smaller for the essential set. While this correlates with the ratio $\frac{\left|\eshpathst\right|}{\left|\shpathst\right|}$ of the size of the essential set to that of the complete set, this correlation is by no means linear: The $82\%$ essential paths out of the complete set of $\Wup{12}$ take approximately $68\%$ of the time of the complete set to compute, while for the $83\%$ of $\Wup{11}$ it takes $76\%$ of the time. 
This phenomenon is even more evident in the two instances \texttt{regional} and \texttt{metro}, in which the complete set takes approximately $20$ and $50$ times of the run time to compute the essential set, respectively, while runtimes are the same for both sets whenever both have the same size. 

However, when looking at individual  computation times, it becomes evident that there are large differences  when considering different source nodes on the same instance: Despite running for approximately 15 minutes on average, for some source nodes the algorithm terminates extremely fast, e.g., in under a second, while a different source node takes more than $127$~minutes for the complete version of \Wup{12}.

Focusing on the number of paths computed, it turns out that the essential set of paths consists of more than $78\%$ of the paths in the complete set for the Wuppertal instances (cf.\ \Cref{tab:all_tests_average}), which indicates that a large majority of the shortest paths are the unique shortest path for some scenario. Moreover, the ratio seems to decrease as the instance size increases. 
The same cannot be observed for the \timpasslib-instances: Here, the two fairly small instances, \texttt{metro} and \texttt{regional} have the lowest essential-to-complete ratio of all instances. 

Intuitively it is obvious that many of the paths in the graph are not shortest paths: Independent of a chosen timetable, some paths via some far away station should never be shortest paths. The extent of this disparity has -- to our knowledge -- not been studied before. 
We therefore compute the total number of paths -- also restricting to paths which are simple when projected onto the $\ptn$ -- and compare this number to our computed complete and essential set of shortest paths. The average number of all paths per source node and per OD-pair can be found in the first two columns of \Cref{tab:all_tests_average}, while the ratio for both the complete and essential set can be found in columns 7 and 8.  

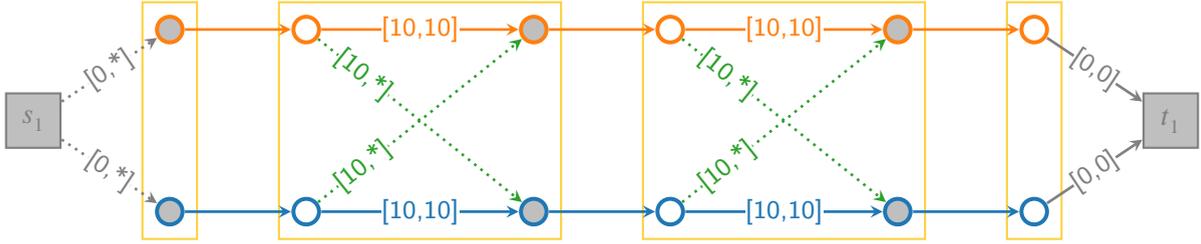
\begin{figure}
    \centering\begin{tikzpicture}[scale = 1.2, transform shape]
	
	\tikzstyle{cell} = [draw, rectangle, gray, thick, fill = gray!50, minimum width = 17, minimum height = 17 ]
	\tikzstyle{a} = [line width=1, -stealth]
	\tikzstyle{transfer} = [a, dotted,tabgreen]
	\tikzstyle{exit} = [a, gray]
	\tikzstyle{enter} = [a, dotted]
	\tikzstyle{l1} = [tabblue]
	\tikzstyle{l2} = [taborange]
	\tikzstyle{l3} = [tabgreen]
	\tikzstyle{border} = [a, line width=0.8, black!60, dotted]
	\tikzstyle{train} = [draw, circle, inner sep=1, minimum width=8 , scale =1, line width = 1.5]
	\tikzstyle{t1} = [train, l1]
	\tikzstyle{t2} = [ train, l2]
	\tikzstyle{t3} = [ train, l3]
	\tikzstyle{path} = [line width=4, zibLightblue]
	\tikzstyle{dep} = [train, black, fill=gray!50]
	\tikzstyle{arr} = [train, black,fill = white]
	\tikzstyle{dep1} = [t1, fill=gray!50]
	\tikzstyle{arr1} = [t1, fill=white]
	\tikzstyle{dep2} = [t2, fill=gray!50]
	\tikzstyle{arr2} = [t2, fill=white]
	\tikzstyle{dep3} = [t3, fill=gray!50]
	\tikzstyle{arr3} = [t3, fill=white]

	\tikzstyle{rop} = [opacity = 0.25]
	
	\tikzstyle{station} = [tabyellow, opacity = .8, thick  ]

    \tikzstyle{midlabel} = [ inner sep = 0, font = \footnotesize, sloped]
    \tikzstyle{me} = [midlabel, gray, fill = white]
    \tikzstyle{mt} = [midlabel, tabgreen, fill = white]
    \tikzstyle{ml1} = [midlabel, l1, fill = white]
    \tikzstyle{ml2} = [midlabel, l2, fill = white]
    \tikzstyle{ml3} = [midlabel, l3, fill = white]
    \tikzstyle{mlt} = [midlabel, gray, fill = white, sloped]

    \draw[station] (1.2,-.3) rectangle (4.3,2.3);
	\draw[station] (-.3,-.3) rectangle (.3,2.3);
	\draw[station] (5.2,-.3) rectangle (8.3,2.3); 
    \draw[station] (9.2,-.3) rectangle (9.8,2.3);

	\node[cell] (s1) at (-1.5,1) {$s_1$};
	\node[cell] (t1) at (11,1) {$t_1$};

	\node[dep1] (AD1) at (0,0)  {};
	\node[arr1] (BA1) at (1.5,0)  {};
	\node[dep1] (BD1) at (4,0) {};
	\node[arr1] (CA1) at (5.5, 0) {};
 	\node[dep1] (CD1) at (8,0) {};
	\node[arr1] (DA1) at (9.5, 0) {};

	\node[dep2] (AD2) at (0,2)  {};
	\node[arr2] (BA2) at (1.5,2)  {};
	\node[dep2] (BD2) at (4,2) {};
	\node[arr2] (CA2) at (5.5, 2) {};
 	\node[dep2] (CD2) at (8,2) {};
	\node[arr2] (DA2) at (9.5, 2) {};

 	\draw[a, l1] (AD1) to (BA1);
    \draw[a, l1] (BA1) to node[ml1] {[10,10]} (BD1);
    \draw[a, l1] (BD1) to (CA1);
    \draw[a, l1] (CA1) to node[ml1] {[10,10]} (CD1);
    \draw[a, l1] (CD1) to (DA1);
    
    \draw[a, l2] (AD2) to (BA2);
    \draw[a, l2] (BA2) to node[ml2] {[10,10]} (BD2);
    \draw[a, l2] (BD2) to (CA2);
    \draw[a, l2] (CA2) to node[ml2] {[10,10]} (CD2);
    \draw[a, l2] (CD2) to (DA2);

    \draw[transfer] (BA1) to node[mt, pos=0.35, left] {[10,*]} (BD2); 
    \draw[transfer] (BA2) to node[mt, pos=0.35, left] {[10,*]} (BD1);
    \draw[transfer] (CA1) to node[mt, pos=0.35, left] {[10,*]} (CD2); 
    \draw[transfer] (CA2) to node[mt, pos=0.35, left] {[10,*]} (CD1);

    \draw[enter] (s1) to node[me] {[0,*]} (AD1);
    \draw[enter] (s1) to node[me] {[0,*]} (AD2);
    \draw[exit] (DA1) to node[me] {[0,0]} (t1);
    \draw[exit] (DA2) to node[me] {[0,0]} (t1);

\end{tikzpicture}
    \caption{Example of a substructure resulting in large differences between the complete and essential set.}
    \label{fig:complete_explosion}
\end{figure}

We observe that in the small instances, \Wup{1}, \Wup{2}, \texttt{toy}, in which almost no routing options are available, (close to) all possible paths are also shortest. However, as soon as the the complexity of the $\ptn$ increases and more transfers occur, only a very small fraction of the paths are indeed shortest, as we see already for \Wup{6} with e.g., $10\%$ of all paths being shortest, while for \Wup{7}, it drops to $2\%$ already. 
From the Wuppertal instances we conclude that, indeed, both the complete and the essential set of shortest paths are significantly smaller than the entire set of paths for realistic instances. 
When considering the \timpasslib{} instances however, the trend can be observed only to some degree: While the fraction of the complete set vs.\ all paths is fairly low for most, the absolute size of the complete set for \texttt{regional} is surprisingly large. Moreover, for the instance \texttt{metro} almost $60\%$ of all paths are in the complete set of shortest paths. At first glance, both numbers seem excessive, particularly when considering also the average values per OD-pair. 
However, these numbers are explainable when taking a closer look at the instances: All dwell times are set to $10$, which coincides with the lower bound on the transfer arcs for \texttt{metro}, while both are set to $1$ in \texttt{regional}. 
Moreover, both contain lines which run in parallel when projected onto \ptn. 
More precisely, both instances contain substructures as in \Cref{fig:complete_explosion}. Suppose both $s_1$-$t_1$-paths which use exclusively the orange or the blue line are shortest. 
Then the essential set will contain only these two paths, yet due to the choice of setting the lower bound on the transfer arcs to the same value as the dwell time, the complete set will contain \emph{each} $s_1$-$t_1$ path in this example. 
As both instances contain multiple such substructures, the cardinality of the essential set differs greatly from that of the complete set, and a large number of all paths are shortest. 
Consequently, this phenomenon can be attributed to the structure of these two specific instances.
The question thus arises whether this modelling of the bounds is realistic. It could be argued that for most modelling purposes dwell times in stations should be shorter than transfer times, making this bound choice unreasonable.

While the computations presented so far have been quite successful, it is also necessary to mention that our algorithm did not terminate for any of our tested source nodes for the instance \texttt{Erding-NDP-S020} within 24 hours. This was rather unexpected, as its size metrics are comparable to some of the Wuppertal instances, which we were able to compute without issues. A closer analysis of \texttt{Erding-NDP-S020} revealed a few disadvantages in the instance structure:  There are multiple lines that run partially in parallel, diverge and reconvene. Moreover, there are express lines, which skip stops, and there is a mix of both high- and low-frequency lines going over the same stops. The lower bounds on the transfer activities are $l = 3$, and upper bounds $u \in \{12, 17, 22, 32, 62\}$ corresponding to lines of frequency $\{6, 4, 3, 2, 1\}$. This combination of factors causes the number of shortest paths to explode: More than $75\%$ of the driving arcs have bounds above $3$, so that a transfer can be favourable, particularly if it goes to a line which skips the next stop. In \texttt{Erding-NDP-S020} the (intended) purpose of our modelling approach backfires: If there exists a fast connection but it runs only once per period, there is a huge amount of alternative paths, which can be options when a vehicle has just left. Thus, what we described in \Cref{ex:frequency_transfer} happens to an extreme in this instance, and too much of a challenge for the algorithm. This interpretation seems to be confirmed by the following small additional test: We manipulated the instance such that all upper bounds on the transfer arcs were at most $32$, i.e., we interpret each line as running at least twice per period time. With this adjustment, we were able to compute both the essential and complete set of shortest paths. We decided against recording these manipulated test results in the table, as we wanted to stick as close to the original instances as possible. 
However, this instance nicely shows the limits of our algorithm -- it is not only dependent on the absolute numbers of lines, arcs, nodes, etc., but also on the specific instance structure and relative bound values. 

\begin{figure}
\vspace{1cm}
\centering
   \begin{subfigure}[t]{.47\textwidth}
       \includegraphics[height=45mm,trim={0 0mm 0 15mm}]
       {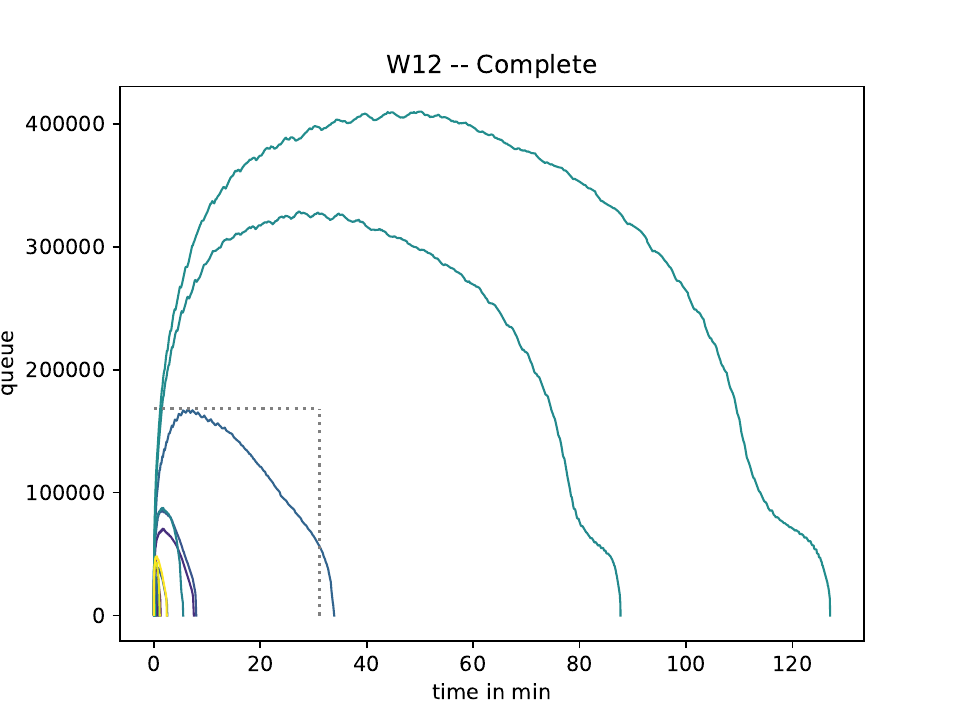}
       \caption{Queue size over time for different source nodes}
       \label{fig:queue_w12_all}
   \end{subfigure}\hfill
   \begin{subfigure}[t]{.47\textwidth}
        \includegraphics[height=45mm,trim={0 1mm 0 9mm}]{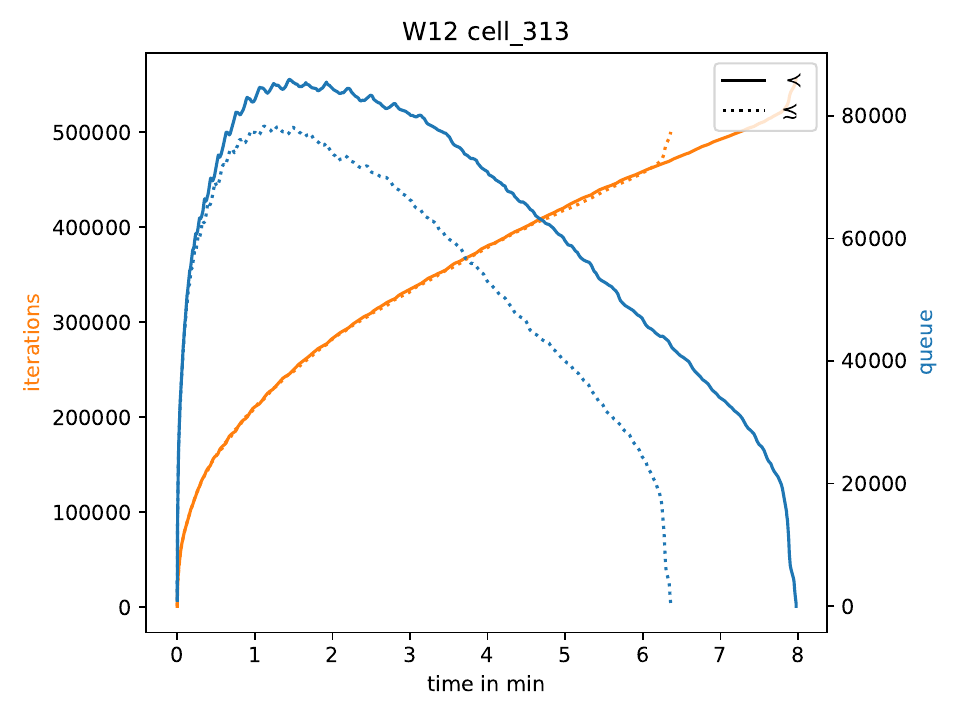}
        \caption{Number of iterations and labels in the queue over time for one source node}
        \label{fig:weak_vs_strict_it_q}
   \end{subfigure}
   \caption{Progress over time of the tropical Dijkstra algorithm}
\end{figure}

In \Cref{fig:queue_w12_all}, we present the queue size over time for \Wup{12} for each source node. All but three runs finish faster than the average time, but the general shapes -- here with one large peak and one smaller bump -- seem to resemble each other.
This can also be observed when comparing both the queue size and iteration number for a complete vs.~essential run of tropical Dijkstra from the same source node, see \Cref{fig:weak_vs_strict_it_q}: At the beginning, both queue and iterations of the complete and essential run coincide, but at some point the complete version seems to shadow the essential run. This behaviour is to be expected, as at the beginning, when not much is explored yet, many tropical path monomials are collected in the labels and inserted into the queue in both cases. In the essential version more monomials can be discarded, which we see as the curves separate.
The complete version does however seem to shadow the behaviour of the essential one, which suggests that the order in which the graph is traversed is broadly similar. 
What \Cref{fig:weak_vs_strict_it_q} also displays nicely is the progression of iterations over time: At the beginning we have a steep incline, as each iteration is easy to do. As the number of dominance checks per iteration increases, the curve flattens slightly, before increasing again at the end, when only few new shortest paths are discovered. 

We conclude that first and foremost, our method is indeed able to compute shortest passenger paths, both the essential set and the complete set -- for decently sized realistic instances, where the enumeration of all paths fails.
Moreover, for realistic instances, the essential set of paths is smaller than the complete set, but not drastically so. This is also reflected by an advantage of computation time.

\subsubsection{Realistic Shortest Passenger Routes}
Our second row of tests is aimed more for the practical application: Usually, passengers perceive transfers as annoying, and journey with large amounts of transfers are undesired. In the previous tests, we did not restrict the number of transfers, and accepted any path which was shortest in the sense of travel time for some scenario. In this row of tests, we want to examine how our algorithm performs when restricting to paths with at most $3$ transfers, as was described in \Cref{subsec:realistic_adjustments}.  

\begin{figure}
    \begin{subfigure}[b]{.48\textwidth}
        \includegraphics[height=45mm]{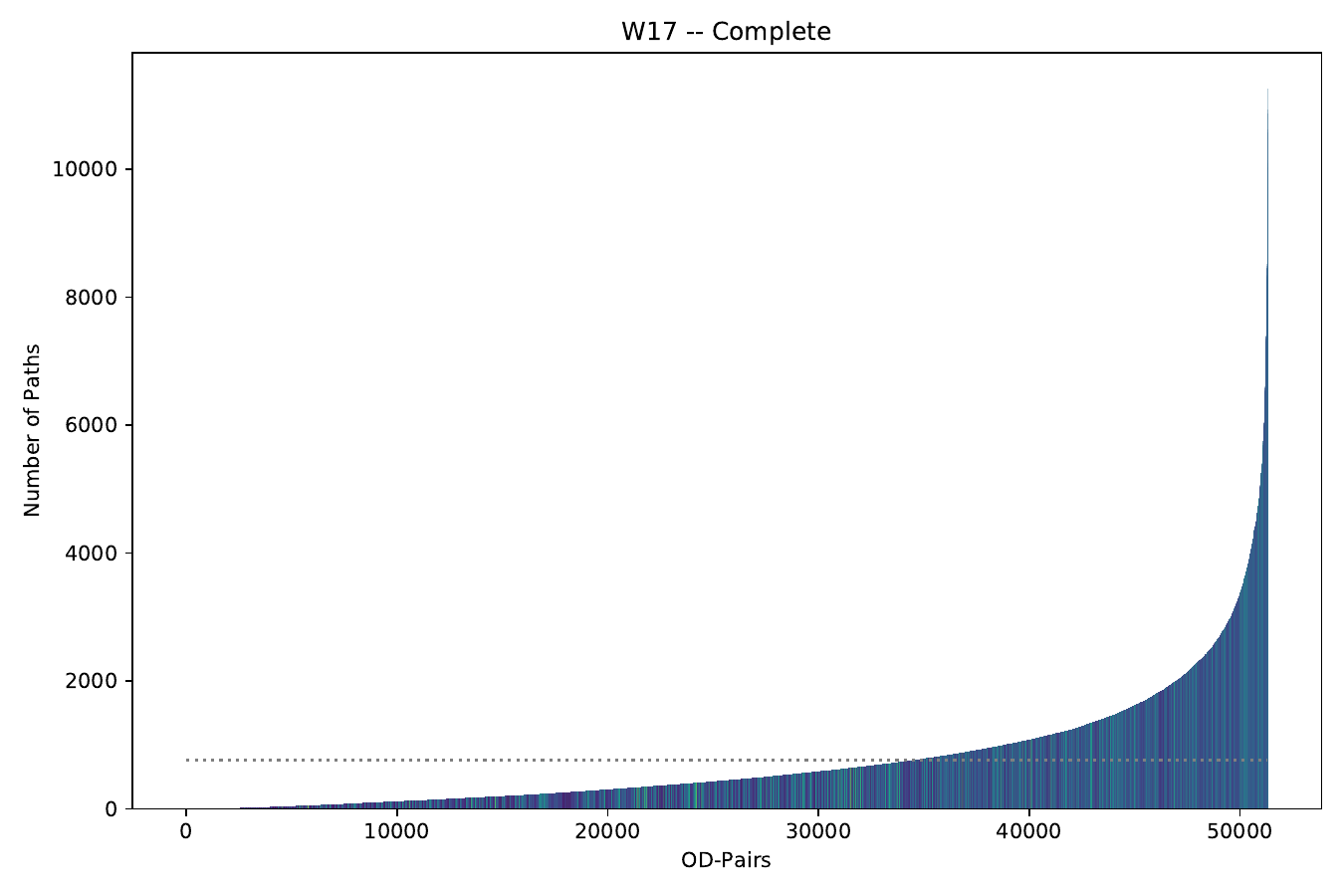}
        \caption{Paths on the event activity network}
    \end{subfigure}
    \hfill
        \begin{subfigure}[b]{.48\textwidth}
        \includegraphics[height=45mm, trim={5mm 3mm 0 10mm}]{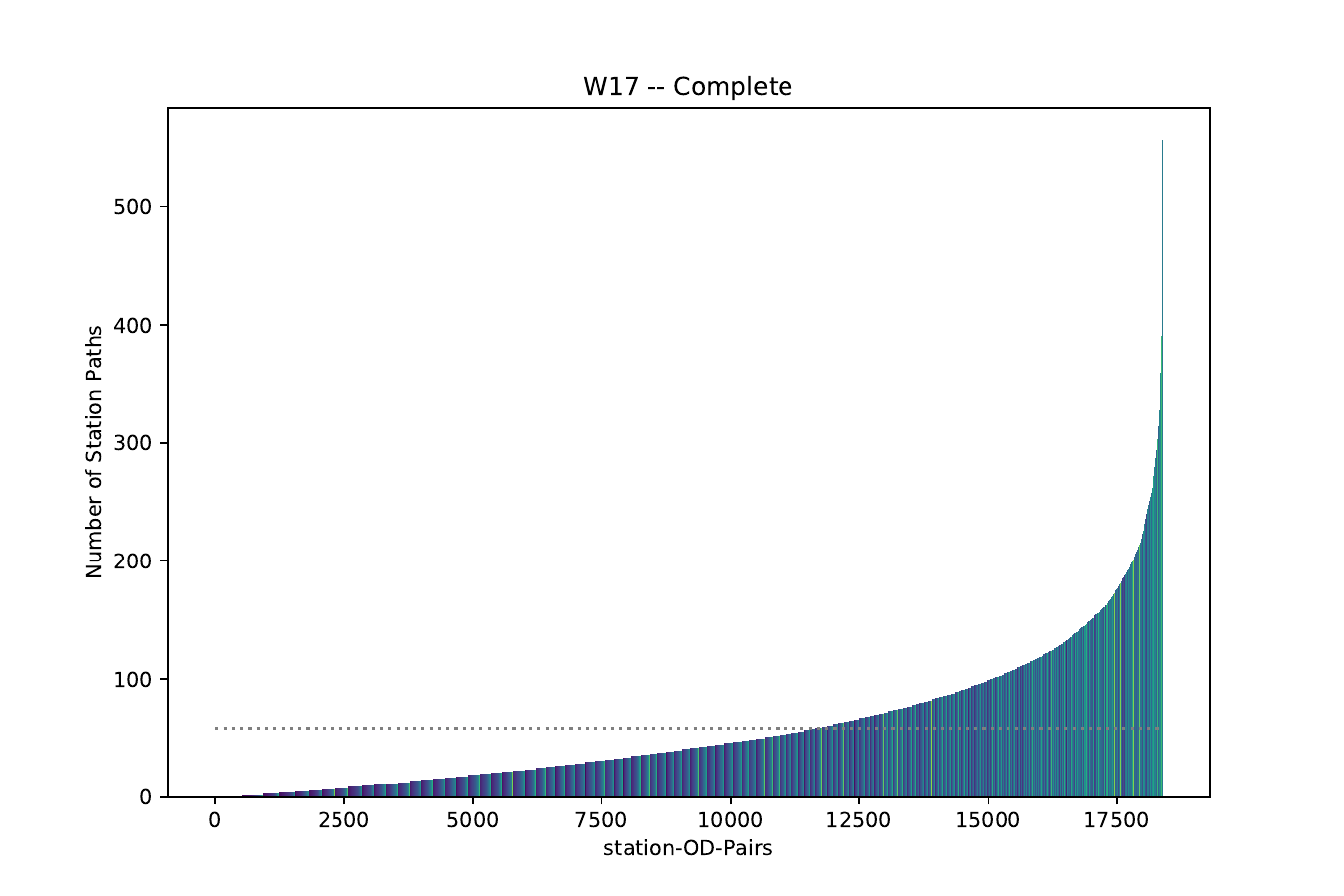}
        \caption{Projection onto the transportation network}
    \end{subfigure}
    \caption{Distribution of number of shortest paths for all OD pairs.}
    \label{fig:distribution_plots}
\end{figure}

With this restriction, we are able to compute both the complete and essential set of paths on a larger row of tests: We computed them for every source node in each of the instances \Wup{1}-\Wup{17}, as well as the first $20$ instances of \timpasslib. For the $6$ largest instances of \timpasslib, we performed the tests only for $20\%$ of the source cells, simply for resource reasons. 
An overview of the results can be found in \Cref{tab:max3wup} and \Cref{tab:max3timpasslib}.
We see that the set of essential paths is always smaller than the complete set, but similar in magnitude. Moreover, the difference in size of the essential and the complete set is not as large as when allowing paths with unrestricted number of transfers. 
For the Wuppertal instances, the number of paths increases with the instance size, but we manage to compute both sets in a very short amount of time: Even for the largest Wuppertal instance \Wup{17} comprising the full transportation network of a city with more than 350,000 inhabitants, we can compute all shortest paths for a source node in less than a minute on average. 
At first glance, the result of over $700$ shortest paths per OD pair on average seems excessive. When taking into account that the upper bound on the transfer times are often very large, cells being connected to multiple stations, and multiple lines running in parallel, this amount is fairly reasonable. 
These aspects motivate to also evaluate how many paths we obtain on the level of the augmented transportation network: For each path in the event-activity network obtained by our algorithm, we project it to its sequence of stations.
The average number of station paths per OD pair of stations can be found in the columns SP/S-od in \Cref{tab:max3wup} and \Cref{tab:max3timpasslib}, which drops down to less than 60 for all instances, even in the complete set. 
However, this average value (both P/st and SP/S-od) is slightly misleading: There are some pairs with very few, and some with many shortest paths. This is evidenced by \Cref{fig:distribution_plots}, which shows the distribution of the number of paths per OD pair.

\begin{table}
    \centering
    
\setlength\tabcolsep{3 pt}
    \begin{tabular}{l | r r r |r |r r r| r |r }
& \multicolumn{4}{c|}{Complete ($\shpathst, \prec$)} & \multicolumn{4}{c|}{Essential ($\eshpathst, \precsim$)}& \\
inst & P/s & P/st & Time/s & SP/S-od & P/s & P/st & Time/s & SP/S-od  &$\%\mathcal C$ \\ \hline
 \Wup{1}  &  $ 23.36 $  &  $ 1.07 $  & $ 0.00 $ ms &  $ 1.00 $  &  $ 23.36 $  &  $ 1.07 $  & $ 0.00 $ ms &  $ 1.00 $  &  $ 100  $  \\ 
 \Wup{2}  &  $ 48.47 $  &  $ 1.61 $  & $ 0.00 $ ms &  $ 1.00 $  &  $ 47.08 $  &  $ 1.57 $  & $ 0.00 $ ms &  $ 1.00 $  &  $ 100  $  \\ 
 \Wup{3}  &  $ 97.60 $  &  $ 2.16 $  & $ 0.00 $ ms &  $ 1.00 $  &  $ 96.90 $  &  $ 2.15 $  & $ 0.00 $ ms &  $ 1.00 $  &  $ 100  $  \\ 
 \Wup{4}  &  $ 158.96 $  &  $ 2.16 $  & $ 0.03 $ ms &  $ 1.00 $  &  $ 156.73 $  &  $ 2.13 $  & $ 0.03 $ ms &  $ 1.00 $  &  $ 100  $  \\ 
 \Wup{5}  &  $ 524.13 $  &  $ 6.76 $  & $ 2.19 $ ms &  $ 1.45 $  &  $ 507.39 $  &  $ 6.54 $  & $ 2.17 $ ms &  $ 1.45 $  &  $ 100  $  \\ 
 \Wup{6}  &  $ 791.84 $  &  $ 8.26 $  & $ 4.95 $ ms &  $ 2.21 $  &  $ 745.81 $  &  $ 7.78 $  & $ 4.74 $ ms &  $ 2.16 $  &  $ 100  $  \\ 
 \Wup{7}  &  $ 1652.28 $  &  $ 15.72 $  & $ 14.20 $ ms &  $ 3.77 $  &  $ 1531.83 $  &  $ 14.57 $  & $ 13.41 $ ms &  $ 3.65 $  &  $ 100  $  \\ 
 \Wup{8}  &  $ 3949.53 $  &  $ 32.16 $  & $ 48.88 $ ms &  $ 7.44 $  &  $ 3711.32 $  &  $ 30.22 $  & $ 46.71 $ ms &  $ 7.35 $  &  $ 100  $  \\ 
 \Wup{9}  &  $ 5604.57 $  &  $ 41.85 $  & $ 82.67 $ ms &  $ 9.00 $  &  $ 5154.79 $  &  $ 38.49 $  & $ 76.71 $ ms &  $ 8.72 $  &  $ 100  $  \\ 
 \Wup{10}  &  $ 9897.83 $  &  $ 67.86 $  & $ 244.76 $ ms &  $ 11.13 $  &  $ 9151.61 $  &  $ 62.75 $  & $ 226.53 $ ms &  $ 10.68 $  &  $ 100  $  \\ 
 \Wup{11}  &  $ 1.31\exn{4} $  &  $ 83.29 $  & $ 387.66 $ ms &  $ 12.20 $  &  $ 1.21\exn{4} $  &  $ 76.90 $  & $ 354.71 $ ms &  $ 11.60 $  &  $ 100  $  \\ 
 \Wup{12}  &  $ 1.78\exn{4} $  &  $ 103.34 $  & $ 678.33 $ ms &  $ 13.32 $  &  $ 1.65\exn{4} $  &  $ 95.84 $  & $ 624.29 $ ms &  $ 12.67 $  &  $ 100  $  \\ 
 \Wup{13}  &  $ 2.79\exn{4} $  &  $ 143.65 $  & $ 1.56 $ s &  $ 20.48 $  &  $ 2.59\exn{4} $  &  $ 133.72 $  & $ 1.45 $ s &  $ 19.57 $  &  $ 100  $  \\ 
 \Wup{14}  &  $ 3.07\exn{4} $  &  $ 147.21 $  & $ 1.95 $ s &  $ 19.09 $  &  $ 2.86\exn{4} $  &  $ 137.14 $  & $ 1.80 $ s &  $ 18.31 $  &  $ 100  $  \\ 
 \Wup{15}  &  $ 4.28\exn{4} $  &  $ 195.81 $  & $ 3.71 $ s &  $ 21.20 $  &  $ 4.00\exn{4} $  &  $ 182.69 $  & $ 3.43 $ s &  $ 20.40 $  &  $ 100  $  \\ 
 \Wup{16}  &  $ 6.61\exn{4} $  &  $ 296.61 $  & $ 9.08 $ s &  $ 32.79 $  &  $ 6.16\exn{4} $  &  $ 276.56 $  & $ 8.44 $ s &  $ 31.64 $  &  $ 100  $  \\ 
 \Wup{17}  &  $ 1.70\exn{5} $  &  $ 763.74 $  & $ 58.40 $ s &  $ 58.31 $  &  $ 1.59\exn{5} $  &  $ 714.96 $  & $ 53.85 $ s &  $ 56.41 $  &  $ 100  $  \\ 
\end{tabular}

    \caption{Average results for each of instances $\Wup{1}-\Wup{17}$ restricted to at most $3$ transfers.  P/s and P/st indicate the average number of paths found per source node and per OD pair, respectively. Time/s is the average time completion of tropical Dijkstra per source node. SP/S-od is the number of station paths per pair of station OD pairs, while $\% \mathcal{C}$ indicates for which percentage of starting cells the tests were performed. }
    \label{tab:max3wup}
\end{table}

When looking at the results of \timpasslib, we do not see the same patterns to the same degree. While the total number of paths in the whole instance increases with instance size, this can be attributed more to the fact that there are also more target cells. When considering the number of paths per OD-pair, this increase with instance size cannot be observed as clearly anymore. 
For example, \texttt{long-distance} is not the largest instance, but has by far the longest average computation time and number of shortest paths, both in the essential and complete sets. Yet the projection onto the \ptn{} reveals only approximately $26$ paths per OD-pair on a station level. 
This is due to the fact that there are many lines running in parallel, in combination with long transfer times, yet the size of the \ptn{} is rather small in comparison to other instances. 
Surprisingly, we were able to compute the instance \texttt{Erding-NDP-S020}, which did not terminate in the previous row of tests, for each starting cell in well under a second on average. 
Here again, we see the relatively large number of shortest paths per OD-pair in comparison to the instance size. The data suggests that the number of paths -- and thus also the computation time -- is larger the more parallel and interconnecting lines there are. This is not particularly surprising, as along parallel lines, the travel times will be similar, and depending on which one the passenger manages to reach on time, more options arise.

\begin{table}
    \centering
    
\setlength\tabcolsep{3 pt}
    \begin{tabular}{l | r r r |r |r r r| r | r}
& \multicolumn{4}{c|}{Complete ($\shpathst, \prec$)} & \multicolumn{4}{c|}{Essential ($\eshpathst, \precsim$)} & \\
inst & P/s & P/st & Time/s   & SP/S-od & P/s & P/st & Time/s & SP/S-od & $\% \mathcal C$ \\ \hline
 \texttt{toy} &  $ 11.63 $  &  $ 1.66 $  & $ 0.00 $ ms &  $ 1.70 $  &  $ 11.50 $  &  $ 1.64 $  & $ 0.00 $ ms &  $ 1.70 $  &  $ 100
 $  \\ 
 \texttt{toy-2} &  $ 69.00 $  &  $ 9.86 $  & $ 2.50 $ ms &  $ 1.70 $  &  $ 68.50 $  &  $ 9.79 $  & $ 2.50 $ ms &  $ 1.70 $  &  $ 100
 $  \\ 
 \texttt{metro} &  $ 1184.20 $  &  $ 23.68 $  & $ 18.59 $ ms &  $ 1.76 $  &  $ 327.22 $  &  $ 6.54 $  & $ 5.16 $ ms &  $ 1.76 $  &  $ 100
 $  \\ 
 \texttt{regional} &  $ 2188.48 $  &  $ 84.17 $  & $ 67.70 $ ms &  $ 2.71 $  &  $ 1209.74 $  &  $ 46.53 $  & $ 49.70 $ ms &  $ 2.71 $  &  $ 100
 $  \\ 
 \texttt{grid} &  $ 1166.64 $  &  $ 48.61 $  & $ 97.16 $ ms &  $ 1.46 $  &  $ 1166.64 $  &  $ 48.61 $  & $ 98.36 $ ms &  $ 1.46 $  &  $ 100
 $  \\ 
 \texttt{Hamburg} &  $ 1001.94 $  &  $ 15.18 $  & $ 23.81 $ ms &  $ 1.79 $  &  $ 958.36 $  &  $ 14.52 $  & $ 22.28 $ ms &  $ 1.77 $  &  $ 100
 $  \\ 
 \texttt{Erding-NDP-S020} &  $ 6424.29 $  &  $ 237.94 $  & $ 699.64 $ ms &  $ 23.76 $  &  $ 6287.96 $  &  $ 232.89 $  & $ 696.46 $ ms &  $ 23.07 $  &  $ 100
 $  \\ 
 \texttt{Schweiz-Fernverkehr} &  $ 6.34\exn{4} $  &  $ 469.41 $  & $ 27.69 $ s &  $ 59.84 $  &  $ 6.15\exn{4} $  &  $ 455.24 $  & $ 26.87 $ s &  $ 59.05 $  &  $ 100
 $  \\ 
 \texttt{R1L1} &  $ 2.39\exn{4} $  &  $ 95.98 $  & $ 1.36 $ s &  $ 16.26 $  &  $ 2.35\exn{4} $  &  $ 94.33 $  & $ 1.34 $ s &  $ 15.95 $  &  $ 100
 $  \\ 
 \texttt{R1L2} &  $ 2.80\exn{4} $  &  $ 112.41 $  & $ 1.64 $ s &  $ 16.69 $  &  $ 2.74\exn{4} $  &  $ 110.15 $  & $ 1.60 $ s &  $ 16.40 $  &  $ 100
 $  \\ 
 \texttt{R1L3} &  $ 2.20\exn{4} $  &  $ 88.57 $  & $ 1.15 $ s &  $ 15.51 $  &  $ 2.17\exn{4} $  &  $ 87.26 $  & $ 1.14 $ s &  $ 15.25 $  &  $ 100
 $  \\ 
 \texttt{R2L1} &  $ 3.63\exn{4} $  &  $ 130.11 $  & $ 2.68 $ s &  $ 19.77 $  &  $ 3.61\exn{4} $  &  $ 129.37 $  & $ 2.67 $ s &  $ 19.64 $  &  $ 100
 $  \\ 
 \texttt{R2L2} &  $ 3.95\exn{4} $  &  $ 141.74 $  & $ 3.09 $ s &  $ 21.45 $  &  $ 3.93\exn{4} $  &  $ 141.01 $  & $ 3.08 $ s &  $ 21.30 $  &  $ 100
 $  \\ 
 \texttt{R1L4} &  $ 3.85\exn{4} $  &  $ 154.48 $  & $ 2.79 $ s &  $ 16.70 $  &  $ 3.79\exn{4} $  &  $ 152.21 $  & $ 2.75 $ s &  $ 16.45 $  &  $ 100
 $  \\ 
 \texttt{R2L3} &  $ 2.91\exn{4} $  &  $ 104.34 $  & $ 1.86 $ s &  $ 18.35 $  &  $ 2.90\exn{4} $  &  $ 103.89 $  & $ 1.85 $ s &  $ 18.25 $  &  $ 100
 $  \\ 
 \texttt{R3L1} &  $ 7.01\exn{4} $  &  $ 237.71 $  & $ 8.05 $ s &  $ 29.51 $  &  $ 6.91\exn{4} $  &  $ 234.20 $  & $ 7.92 $ s &  $ 29.15 $  &  $ 100
 $  \\ 
 \texttt{R3L2} &  $ 7.65\exn{4} $  &  $ 259.20 $  & $ 9.79 $ s &  $ 31.30 $  &  $ 7.53\exn{4} $  &  $ 255.25 $  & $ 9.66 $ s &  $ 30.92 $  &  $ 100
 $  \\ 
 \texttt{R4L1} &  $ 8.54\exn{4} $  &  $ 268.45 $  & $ 10.55 $ s &  $ 43.59 $  &  $ 8.42\exn{4} $  &  $ 264.84 $  & $ 10.40 $ s &  $ 43.04 $  &  $ 100
 $  \\ 
 \texttt{R4L2} &  $ 9.37\exn{4} $  &  $ 294.78 $  & $ 12.54 $ s &  $ 42.10 $  &  $ 9.26\exn{4} $  &  $ 291.27 $  & $ 12.42 $ s &  $ 41.53 $  &  $ 100
 $  \\ 
 \texttt{R3L3} &  $ 6.91\exn{4} $  &  $ 234.36 $  & $ 7.65 $ s &  $ 27.80 $  &  $ 6.80\exn{4} $  &  $ 230.54 $  & $ 7.53 $ s &  $ 27.52 $  &  $ 100
 $  \\ 
 \texttt{long-distance} &  $ 7.13\exn{5} $  &  $ 2985.35 $  & $ 3088.12 $ s &  $ 26.93 $  &  $ 7.02\exn{5} $  &  $ 2938.60 $  & $ 3001.61 $ s &  $ 26.35 $  &  $ 20
 $  \\ 
 \texttt{Stuttgart} &  $ 4.10\exn{5} $  &  $ 1057.37 $  & $ 1114.01 $ s &  $ 38.22 $  &  $ 4.09\exn{5} $  &  $ 1055.92 $  & $ 1143.18 $ s &  $ 38.19 $  &  $ 20
 $  \\ 
 \texttt{R4L3} &  $ 1.04\exn{5} $  &  $ 328.01 $  & $ 14.71 $ s &  $ 44.53 $  &  $ 1.03\exn{5} $  &  $ 324.29 $  & $ 14.43 $ s &  $ 44.01 $  &  $ 20
 $  \\ 
 \texttt{R2L4} &  $ 7.88\exn{4} $  &  $ 282.54 $  & $ 9.46 $ s &  $ 21.19 $  &  $ 7.83\exn{4} $  &  $ 280.68 $  & $ 9.47 $ s &  $ 21.06 $  &  $ 20
 $  \\ 
 \texttt{R3L4} &  $ 1.32\exn{5} $  &  $ 446.14 $  & $ 24.29 $ s &  $ 30.40 $  &  $ 1.30\exn{5} $  &  $ 440.76 $  & $ 23.92 $ s &  $ 30.07 $  &  $ 20
 $  \\ 
 \texttt{R4L4} &  $ 1.83\exn{5} $  &  $ 576.65 $  & $ 36.56 $ s &  $ 49.22 $  &  $ 1.80\exn{5} $  &  $ 567.46 $  & $ 35.87 $ s &  $ 48.55 $  &  $ 20
 $  \\ 
\end{tabular}

    \caption{Average results for each of the \timpasslib instances, restricted to at most $3$ transfers.  P/s and P/st indicate the average number of paths found per source node and per OD pair, respectively. Time/s is the average time completion of tropical Dijkstra per source node. SP/S-od is the number of station paths per pair of station OD pairs, while $\% \mathcal{C}$ indicates for which percentage of starting cells the tests were performed.}
    \label{tab:max3timpasslib}
\end{table}

Apart from just evaluating our algorithm, we also wanted to gain some insight into how many paths are potentially relevant when it comes to shortest passenger paths in the application of timetabling.
We have to admit that we were surprised by the number of shortest paths. 
Without restricting the number of transfers, there is a larger difference between the complete, and the essential set, but with the restriction, both sets seem to have a similar size for practical instances. 
Nevertheless, even the essential set of shortest paths is in most cases quite large, even when projected onto the underlying \ptn. 
Our results further support the need to consider timetabling and passenger routing in combination, as a different timetable -- corresponding to a scenario in our case -- can lead to entirely different best passenger paths.
Fixing passenger paths in advance can then lead to overall suboptimal timetables.

\section{Conclusion and Outlook}
\label{sec:conclusion}

The tropical Dijkstra is a powerful technique to solve \cispp and \eispp on graphs with interval costs. From our case study, we conclude that it is possible to compute both the complete and an essential set of shortest paths on realistic transportation networks. 
Restricting the paths to be meaningful for this application has sped up the computation significantly and enables to solve considerably larger instances. 
Our results show the relevance of considering timetabling and passenger routing in an integrated manner. 
Moreover, we hope our procedure will be fruitful to foster algorithms for integrated timetabling and passenger routing in the future. 
For example, computing an essential set as preprocessing can avoid re-computation of shortest paths in iterative modulo network simplex approaches \cite{lobel_restricted_2020}.

\bibliographystyle{plain}

\bibliography{references}

\begin{thebibliography}{10}

\bibitem{catanzaro_reduction_2011}
D.~Catanzaro, M.~Labbé, and M.~Salazar-Neumann.
\newblock Reduction approaches for robust shortest path problems.
\newblock {\em Computers \& Operations Research}, 38(11):1610--1619, 2011.

\bibitem{martins_multicriteria_1984}
E.~de~Queirós Vieira~Martins.
\newblock On a multicriteria shortest path problem.
\newblock {\em European Journal of Operational Research}, 16(2):236--245, 1984.

\bibitem{ehrgott}
M.~Ehrgott.
\newblock {\em Multicriteria Optimization}.
\newblock Springer-Verlag, Berlin, Heidelberg, 2005.

\bibitem{hamonic_optimizing_2023}
F.~Hamonic, C.~Albert, B.~Couëtoux, and Y.~Vaxès.
\newblock Optimizing the ecological connectivity of landscapes.
\newblock {\em Networks}, 81(2):278--293, 2023.

\bibitem{joswig_essentials_trop}
M.~Joswig.
\newblock {\em Essentials of tropical combinatorics}, volume 219 of {\em
  Graduate Studies in Mathematics}.
\newblock American Mathematical Society, Providence, RI, 2021.

\bibitem{JoswigSchroeter_2022}
M.~Joswig and B.~Schröter.
\newblock Parametric shortest-path algorithms via tropical geometry.
\newblock {\em Mathematics of Operations Research}, 47(3):2065–2081, 2022.

\bibitem{karasan_robust_2012}
O.~E. Karasan, M.~Pinar, and H.~Yaman.
\newblock The robust shortest path problem with interval data.
\newblock Technical report, Bilkent University, 2001.

\bibitem{lindner_optimal_2021}
N.~Lindner, P.~Maristany de~las Casas, and P.~Schiewe.
\newblock Optimal {Forks}: Preprocessing single-source shortest path instances
  with interval data.
\newblock In M.~Müller-Hannemann and F.~Perea, editors, {\em 21st {Symposium}
  on {Algorithmic} {Approaches} for {Transportation} {Modelling},
  {Optimization}, and {Systems} ({ATMOS} 2021)}, volume~96 of {\em Open
  {Access} {Series} in {Informatics} ({OASIcs})}, pages 7:1--7:15, Dagstuhl,
  Germany, 2021. Schloss Dagstuhl – Leibniz-Zentrum für Informatik.

\bibitem{lobel_restricted_2020}
F.~Löbel, N.~Lindner, and R.~Borndörfer.
\newblock The restricted modulo network simplex method for integrated periodic
  timetabling and passenger routing.
\newblock In J.~S. Neufeld, U.~Buscher, R.~Lasch, D.~Möst, and
  J.~Schönberger, editors, {\em Operations {Research} {Proceedings} 2019},
  pages 757--763, Cham, 2020. Springer International Publishing.

\bibitem{tropicalGeom_sturmfels}
D.~Maclagan and B.~Sturmfels.
\newblock {\em Introduction to Tropical Geometry}.
\newblock Providence: American Mathematical Society, 2015.

\bibitem{maristany_de_las_casas_new_2024}
P.~Maristany de~las Casas.
\newblock {\em New Multiobjective Shortest Path Algorithms}.
\newblock PhD thesis, Freie Universität Berlin, Berlin, 2024.

\bibitem{maristany_de_las_casas_improved_2021}
P.~Maristany de~las Casas, A.~Sedeño-Noda, and R.~Borndörfer.
\newblock An improved multiobjective shortest path algorithm.
\newblock {\em Computers \& Operations Research}, 135:105424, 2021.

\bibitem{masing_forward}
B.~Masing, N.~Lindner, and P.~Ebert.
\newblock Forward and line-based cycle bases for periodic timetabling.
\newblock {\em Oper. Res. Forum}, 4(3), June 2023.

\bibitem{montemanni_exact_2004}
R.~Montemanni and L.~M. Gambardella.
\newblock An exact algorithm for the robust shortest path problem with interval
  data.
\newblock {\em Computers \& Operations Research}, 31(10):1667--1680, 2004.

\bibitem{schiewe_introducing_2023}
P.~Schiewe, M.~Goerigk, and N.~Lindner.
\newblock Introducing {TimPassLib} – a library for integrated periodic
  timetabling and passenger routing.
\newblock {\em Operations Research Forum}, 4(3):64, 2023.

\bibitem{schiewe_periodic_2020}
P.~Schiewe and A.~Schöbel.
\newblock Periodic timetabling with integrated routing: Toward applicable
  approaches.
\newblock {\em Transportation Science}, 54(6):1714--1731, 2020.
\newblock Publisher: INFORMS.

\bibitem{schmidt_integrating_2014}
M.~Schmidt.
\newblock {\em Integrating Routing Decisions in Public Transportation
  Problems}, volume~89 of {\em Springer {Optimization} and {Its}
  {Applications}}.
\newblock Springer, New York, NY, 2014.

\end{thebibliography}
\newpage

\appendix

\section{Proofs}

\technical*
\begin{proof}
    Let $p, p' \in \pathst$.
    \begin{enumerate}
        \item Suppose that $c(p) \vartriangleleft c(p')$. By construction of the best scenario, we have
        \begin{align}
        \bestreal{p}(p) &= \sum_{a\in \arcs(p)\setminus \arcs(p')} \ell_a + \sum_{a\in \arcs(p)\cap\arcs(p')} \ell_a
        \leq \sum_{a\in \arcs(p)\setminus \arcs(p')} c_a + \sum_{a\in \arcs(p)\cap\arcs(p')} \ell_a \\
        &
        \vartriangleleft \sum_{a\in \arcs(p')\setminus \arcs(p)} c_a + \sum_{a\in \arcs(p)\cap\arcs(p')} \ell_a \leq \sum_{a\in \arcs(p')\setminus \arcs(p)} u_a + \sum_{a\in \arcs(p)\cap\arcs(p')} \ell_a \\
        &= \bestreal{p}(p'),
    \end{align}
    where we use that $c(p) \vartriangleleft c(p')$ implies
    \begin{equation}
        \sum_{a\in \arcs(p)\setminus \arcs(p')} c_a
        = 
        c(p) - \sum_{a\in \arcs(p)\cap \arcs(p')} c_a
        \;\vartriangleleft\;
         c(p') - \sum_{a\in \arcs(p)\cap \arcs(p')} c_a
        = \sum_{a \in \arcs(p')\setminus \arcs(p)} c_a .
    \end{equation}
    \item Let $c \in [\bl, \bu]$ and suppose that $\bestreal{p'}(p) \vartriangleleft \bestreal{p'}(p')$. First, observe 
    \begin{equation}\label{lb}
        \sum_{a\in \arcs(p)\setminus \arcs(p')} c_a \leq \sum_{a\in \arcs(p) \setminus \arcs(p') } u_a = \sum_{a\in \arcs(p)\setminus \arcs(p')} \bestreal{p'}_a, 
    \end{equation}
    while 
    \begin{equation}\label{ub}
        \sum_{a\in \arcs(p') \setminus \arcs(p)} c_a \geq \sum_{a\in \arcs(p') \setminus \arcs(p)} \ell_a 
        = \sum_{a\in \arcs(p') \setminus \arcs(p)} \bestreal{p'}_a.
    \end{equation}
    Moreover, since $\bestreal{p'}(p) \vartriangleleft \bestreal{p'}(p')$ it must hold that
    \begin{equation}\label{cmpbest}
        \sum_{a\in \arcs(p)\setminus \arcs(p')} \bestreal{p'}_a \vartriangleleft \sum_{a\in \arcs(p')\setminus \arcs(p)} \bestreal{p'}_a.
    \end{equation}
    
    Consequently, for the chosen $c$, we obtain 
    \begin{align}
        c(p)    &= \sum_{a\in \arcs(p)\setminus \arcs(p')} c_a + \sum_{a\in \arcs(p) \cap \arcs(p')} c_a \\
        &\overset{\eqref{lb}}{\leq} \sum_{a\in \arcs(p)\setminus \arcs(p')} \bestreal{p'}_a + \sum_{a\in \arcs(p) \cap \arcs(p')} c_a\\
        & \overset{\eqref{cmpbest}}{\vartriangleleft} \sum_{a\in \arcs(p')\setminus \arcs(p)} \bestreal{p'}_a + \sum_{a\in \arcs(p) \cap \arcs(p')} c_a  \\ &\overset{\eqref{ub}}{\leq} \sum_{a\in \arcs(p') \setminus \arcs(p)} c_a + \sum_{a\in \arcs(p) \cap \arcs(p')} c_a \\
        &= c(p'). \qedhere
    \end{align}
    \end{enumerate}
\end{proof}

\constantpaths*
\begin{proof}
    Clearly (1) $\implies$ (2) $\implies$ (3). By \Cref{lem:technical}, (3) implies (1). Now suppose that $p \sim p'$. Then, by (2), 
    \begin{equation}
         \sum_{a\in \arcs(p) \setminus \arcs(p')} \ell_a = \sum_{a\in \arcs(p) \setminus \arcs(p')} \bestreal{p}_a = \sum_{a\in \arcs(p) \setminus \arcs(p')} \bestreal{p'}_a = \sum_{a\in \arcs(p) \setminus \arcs(p')} u_a,
    \end{equation}
    so that $\ell_a = u_a$ for all $a \in \arcs(p) \setminus \arcs(p')$. By an analogous argument, $\ell_a = u_a$ for all $a \in \arcs(p') \setminus \arcs(p)$.
\end{proof}

\equivalencenonshortestweakpath*
\begin{proof}
    (1) $\implies$ (2): Let $p, p'\in \eshpathst$ with $p \neq p'$. If $p' \preceq p$ or $p \sim p'$, then $c(p') \leq c(p)$ for all $c \in [\bl, \bu]$. So $\eshpathst \setminus \{p\}$ is still scenario-covering, and hence $\eshpathst$ cannot be essential, which leads to a contradiction.

    (2) $\implies$ (3): By \Cref{lem:bestreal_weak_dom}, we have $\bestreal{p}(p') > \bestreal{p}(p)$ or $\bestreal{p'}(p') \geq \bestreal{p'}(p)$. Moreover, \Cref{lem:constant_paths} implies $\bestreal{p}(p') > \bestreal{p}(p)$ or $\bestreal{p'}(p) > \bestreal{p'}(p')$, which means that (3) must hold. 
    
    (3) $\implies$ (1): If $p \notin \eshpathst$, then $p$ is by (3) shorter than any other path in $\eshpathst$ w.r.t.\ $\bestreal{p}$, so $\eshpathst$ is not scenario-covering, thus contradicting the assumption.
\end{proof}
\end{document}